\documentclass[a4paper,reqno]{amsart}

\usepackage{amsmath}
\usepackage{amsthm}
\usepackage{amsfonts}
\usepackage{amssymb}

\usepackage{hyperref}

\hypersetup{
    colorlinks=false,
    pdfborder={0 0 0},
    pdfstartview={XYZ null null 1.00}
}

\numberwithin{equation}{section}

  \def\degl{\mathrm{deg_{LS}}}
  
  \def\degb{\mathrm{deg_B}}

  \def\<{\langle}
  \def\>{\rangle}

  \def\ker{\mathrm{Ker}\,}
  \def\im{\mathrm{Im}\,}

  \def\re{\mathrm{Re}\,}
  \def\ve{\varepsilon}

  \newcommand{\y}{\bar v}
  \newcommand{\z}{\bar w}

  \def\d{ \, d }

  \def\R{\mathbb{R}}

  \def\C{\mathbb{C}}

  \def\ffor{ \qquad \mathrm{for} \quad }
  \def\aas{ \qquad \mathrm{as} \quad }

  \def\o{\overline}
  \def\t{\widetilde}
  
  \def\h{\widehat}

\theoremstyle{plain}
  \newtheorem{theorem}{Theorem}[section]

  \newtheorem{proposition}[theorem]{Proposition}
  
  \newtheorem{lemma}[theorem]{Lemma}
  \newtheorem{corollary}[theorem]{Corollary}

\theoremstyle{definition}
  \newtheorem{definition}[theorem]{Definition}
  
  \newtheorem{remark}[theorem]{Remark}

\begin{document}

\title[Averaging principle and periodic solutions...]{Averaging principle and periodic solutions for nonlinear evolution equations at resonance}

\author{Piotr Kokocki}
\address{\noindent Faculty of Mathematics and Computer Science \newline Nicolaus Copernicus University \newline Chopina 12/18, 87-100 Toru\'n, Poland}
\email{pkokocki@mat.umk.pl}
\thanks{The researches supported by the MNISzW Grant no. N N201 395137}



\begin{abstract}
We study the existence of $T$-periodic solutions $(T > 0)$ for the first order differential equations being at resonance at infinity, where the right hand side is the perturbations of a sectorial operator. Our aim is to prove an index formula expressing the topological degree of the associated translation along trajectories operator on appropriately large ball, in terms of special geometrical assumptions imposed on the nonlinearity. We also prove that the geometrical assumptions are generalization of well known Landesman-Lazer and strong resonance conditions.
Obtained index formula is used to derive the criteria determining the existence of $T$-periodic solutions for the heat equation being at resonance at infinity.
\end{abstract}

\maketitle

\setcounter{tocdepth}{2}

\section{Introduction}

We consider nonlinear differential equations of the form
\begin{align}\label{row-par}
\dot u(t) = - A u(t) + \lambda u(t) + F (t, u(t)),  \qquad  t > 0
\end{align}
where $\lambda$ is a real number, $A: X\supset D(A)\to X$ is a sectorial operator on a Banach space $X$ and $F:[0,+\infty)\times X^\alpha\to X$ is a continuous map. Here $X^\alpha$ for $\alpha\in(0,1)$, is a fractional power space given by $X^\alpha:=D((A + \delta I)^\alpha)$, where $\delta > 0$ is such that the operator $A + \delta I$ is  positively defined.
This equation is an abstract formulation of many partial differential equations including the nonlinear heat equation
\begin{equation}\label{eq2wave}
    u_t(x,t) = \Delta u(x,t) + \lambda u(x,t) + f(t,x,u(x,t)) \ffor t\ge 0, \ x\in\Omega
\end{equation}
where $\Omega$ is an open subset of $\R^n$ ($n \ge 1$), $\Delta$ is a Laplace operator with the Dirichlet boundary conditions and $f:[0,+\infty)\times\Omega\times\R\to\R$ is a continuous map. To see this, it is enough to take $A u:= - \Delta u$ and $F(t,u) = f(t,\,\cdot\,, u(\cdot))$.

In this paper, we intend to study the existence of $T$-periodic solutions ($T > 0$) for the equation \eqref{row-par} being at {\em resonance at infinity}, that is, $\ker(\lambda I - A)\neq \{0\}$ and $F$ is a bounded map. To explain this more precisely assume that, for every initial data $x\in X^\alpha$, the equation \eqref{row-par} admits a (mild) solution $u:[0,+\infty)\to X^\alpha$ starting at $x$. Then the $T$-periodic solutions of \eqref{row-par} can be identified with fixed points of \emph{the translation along trajectories operator} $\Phi_T:X^\alpha \to X^\alpha$, defined by
\begin{equation*}
\Phi_T (x):= u(T;x) \ffor x\in X^\alpha
\end{equation*}
Effective methods for studying the existence of fixed points of translation along trajectories operator are so called \emph{averaging principles}, expressing the fixed point index of the operator $\Phi_T$ in terms of the averaging of the right side of \eqref{row-par}. If the topological degree of this averaging is nontrivial, then the translation along trajectories operator admits a fixed point, which in turn is a starting point of $T$-periodic solution.

The averaging principle for equations on finite dimensional manifolds were studied in \cite{MR819192}, while the generalization on the case of equations on arbitrary Banach spaces were considered in \cite{MR2183379}, when the right side of equation is a nonlinear perturbation of a generator of $C_0$ semigroup. In \cite{cw-kok1} the averaging principle were studied in the case when $-A$ generates a $C_0$ semigroup of contractions and $F$ is a condensing map with respect to the Hausdorff measure of noncompactness. The results for equations with the right side being a nonlinear perturbation of the family of generators of $C_0$ semigroups $\{A(t)\}_{t\ge 0}$ are contained in \cite{cwi-kok2}.

{\em The resonant version of averaging principle} was proved in \cite{Kok1} in the case when $A$ is a generator of a compact $C_0$ semigroup (not necessary sectorial) and $F:[0,+\infty)\times X\to X$ is a continuous map. Obtained result were used to prove the criteria on the existence of $T$-periodic solution for \eqref{eq2wave} under the assumption that $f$ satisfies Landesman-Lazer conditions. In this paper we continue the studies from \cite{Kok1}. First we prove the resonant version of averaging principle for the equation \eqref{row-par} in the case when $A$ is a sectorial operator. This assumption is stronger than that in \cite{Kok1}, however it allows us considering the wider class of nonlinear perturbations, namely we permit the maps $F:[0,+\infty)\times X^\alpha \to X$ defined on fractional power spaces with $\alpha\in[0,1)$. Subsequently we use this principle to prove criteria determining the existence of $T$-periodic solutions for \eqref{row-par}, in the terms of appropriate geometrical conditions imposed on $F$. It turns out that these geometrical conditions are more general that Landesman-Lazer conditions used in \cite{Kok1}.

The main difficulty lies in the fact that, in the presence of resonance, the problem of existence of periodic solutions may not have solution for general nonlinearity $F$. This fact has been explained in detail in Remark \ref{rem-non-ex2}. We overcome this difficulty providing new theorems determining the existence of $T$-periodic solutions of \eqref{row-par}, in terms of appropriate geometrical assumptions imposed on the nonlinearity $F$. To formulate this geometrical assumptions we will need a special direct sum decomposition of the space $X:=X_-\oplus X_0\oplus X_+$, which will be obtained in Theorem \ref{th:10} as the main result of Section 2. This direct sum decomposition is actually a spectral decomposition of the operator $A$ with the property that $X_0:=\ker(\lambda I - A)$, and the parts of the operator $\lambda I - A$ in $X_+$ and $X_-$ are positively and negatively defined, respectively.

\noindent Section 3 is devoted to the mild solutions for \eqref{row-par}. First we remind the standard facts concerning the existence and uniqueness for this equation and then we discuss the continuity of mild solutions with respect to the initial data and parameter. Furthermore, as we will use the homotopy invariants, we provide some compactness properties for the translation operator. More precisely, we prove that $\Phi_T$ is completely continuous map provided $A$ has compact resolvents.

\noindent In Section 4 we prove the first result: the resonant version of averaging principle. More precisely, we will consider the equations of the form
\begin{align}\label{row-par-ep}
\dot u(t) = - A u(t) + \lambda u(t) + \ve F (t, u(t)),  \qquad  t > 0
\end{align}
where $\ve\in[0,1]$ is a parameter. Let $\Phi_T(\ve,\,\cdot\,):X^\alpha\to X^\alpha$ be the associated translation along trajectories operator and let $g\colon N_\lambda \to N_\lambda$, where $N_\lambda:= \ker (\lambda I - A)$, be a map given by
\begin{equation*}
g(x) := \int_0^T PF(\tau,x)\,d \tau \qquad\mathrm{for}\quad x\in N_\lambda.
\end{equation*}
Write $X^\alpha_+:= X^\alpha\cap X_+$, $X^\alpha_-:= X^\alpha\cap X_-$ and assume that $U\subset N_\lambda$ and $V\subset X^\alpha_-\oplus X^\alpha_+$ are such that $0\in V$ and $g(x)\neq 0$ for $x\in\partial U$.
Then {\em the resonant averaging principle} says that, for small $\ve > 0$, the fixed point index of $\Phi_T(\ve,\,\cdot\,)$ is equal to the Brouwer degree of $-g$.

\noindent In Section 5 we formulate geometrical conditions $(G1)$ and $(G2)$ (see page \pageref{g1g2}) and apply the resonant averaging principle to prove the second result, {\em the index formula for periodic solutions}, which express the fixed point index of the translation along trajectories operator $\Phi_T$ on sufficiently large ball, in terms of conditions $(G1)$ and $(G2)$.

\noindent Finally, in Section 6 we provide applications for particular partial differential equations. First of all, in Theorems \ref{lem-est2} and \ref{lem-est3}, we prove that if $F$ is a Niemytzki operator associated with a map $f$, then the well known in literature Landesman-Lazer (see e.g.\cite{MR0267269}, \cite{MR0513090}) and strong resonance conditions (see e.g. \cite{MR713209}, \cite{MR597281}, \cite{MR1055536}) are actually particular cases of assumptions $(G1)$ and $(G2)$.

\section{Spectral decomposition of linear operators}

Let $A:X\supset D(A)\to X$ be a sectorial operator on a real Banach space $X$ with a norm $\|\cdot\|$ such that: \\[5pt]
\noindent\makebox[9mm][l]{$(A1)$}\parbox[t][][t]{118mm}{the resolvent of the operator $A$ are compact,}\\[5pt]
\noindent\makebox[9mm][l]{$(A2)$}\parbox[t][][t]{118mm}{there is an injection $i:X \hookrightarrow H$, where $H$ is a Hilbert space with norm $\|\cdot\|_H$ and scalar product $\<\,\cdot\,, \,\cdot\,\>_H$,}\\[5pt]
\noindent\makebox[9mm][l]{$(A3)$}\parbox[t][][t]{118mm}{there exists a self-adjoint operator $\h A:H\supset D(\h A) \to H$ such that $$i\times i\left[\mathrm{Gr}\,(A)\right]\subset \mathrm{Gr}\,(\h A).$$ }\\[5pt] 
  
Since $A:X\supset D(A)\to X$ is a sectorial operator, there is $\delta \ge 0$ such that $\re z > 0$ for $z\in\sigma(A + \delta I)$. Write $A_\delta:= A + \delta I$. Our aim in this section is to prove the following theorem.

\begin{theorem}\label{th:10ab}
Assume that $(A1)$, $(A2)$ and $(A3)$ hold. If $\lambda = \lambda_k$ for some $k\ge 1$, is an eigenvalue of $A$, then there is a direct sum decomposition on closed subspaces $X = X_+\oplus X_-\oplus X_0$ such that
$$S_A(t)X_i \subset X_i \ffor t\ge 0, \ i\in\{0,-,+\},$$ and the following assertions are satisfied:
\begin{enumerate}
\item[(i)] $X_0 = \ker (\lambda I - A)$, $X_-$ is a finite dimensional space such that
\begin{equation*}
   X_-= \{0\} \text{ if } k=1 \  \text{ and } \   X_-=\bigoplus_{i=1}^{k-1} \ker(\lambda_i I - A) \text{ if } k\ge 2.
\end{equation*}
Hence $\dim X_- = 0$ if $k=1$ and $\dim X_- = \sum_{i=1}^{k-1} \dim\ker(\lambda_i I - A)$ if $k\ge 2$.
\item[(ii)] we have the following inequalities 
\begin{align}\label{ine11ab} 
     & \|A_\delta^\alpha S_A(t)x\| \le K e^{- (\lambda + c) t} t^{-\alpha} \|x\| &&     \text{for} \quad x\in X_+, \  t > 0,\\ \label{ine22ab}
     & \|e^{\lambda t}S_A(t)x\| \le K e^{- c t}\|x\| && \text{for} \quad x\in X_+, \ t\ge 0, \\ \label{ine33ab}
     & \|e^{\lambda t}S_A(t)x\| \le K e^{c t}\|x\| &&\text{for} \quad x\in X_-,\ t\le 0,
\end{align}  
where $c, K > 0$ are constants,
\item[(iii)]we have the following orthogonality condition $$\<i (u_l),i(u_m)\>_H = 0$$ for $u_l\in X_l$ and $u_m\in X_m$ where $l,m\in\{0,-,+\}$, $l\neq m$.  
\end{enumerate}
\end{theorem}

Before start the proof we recall that \emph{the complexification} of the linear space $X$ is, by definition, a complex linear space $X_\C:=X\times X$ with the following operations 
\begin{align*}   
& (z_1,z_2) + (z'_1,z'_2) = (z_1 + z'_2, z_1 + z'_2) && \text{ if }\ \  (z_1,z_2), (z'_1,z'_2)\in X_\C, \\
& \lambda\cdot (z_1,z_2) = (\lambda_1 z_1 - \lambda_2 z_2, \lambda_1 z_2 + \lambda_2 z_1) && \text{ if }\ \  \lambda=(\lambda_1 + \lambda_2 i)\in\C, \ (z_1,z_2)\in X_\C.  
\end{align*}
Writing $z_1 + i z_2:=(z_1,z_2)$ for $(z_1,z_2)\in X_\C$, the above operations take the natural form
\begin{equation*}
\begin{aligned}
& (z_1 + i z_2) + (z'_1 + i z'_2) = (z_1 + z'_1) + i(z_2 + z'_2) \\
& \lambda\cdot (z_1 + iz_2) = (\lambda_1 z_1 - \lambda_2 z_2) + i(\lambda_1 z_2 + \lambda_2 z_1)
\end{aligned}
\end{equation*} 
for $z_1 + i z_2, z'_1 + i z'_2\in X_\C$ and $\lambda=(\lambda_1 + \lambda_2 i)\in\C$. We recall that the complexification of an operator $A:X\supset D(A)\to X$ defined on a real linear space $X$ is a linear operator $A_\C: X_\C \supset D(A_\C) \to X_\C$ given by
\begin{align*}
& D(A_\C):= \{z_1 + i z_2 \in X_\C \ | \ z_1, z_2\in D(A)\}, \\
& A_\C z := A z_1 + i A z_2 \ffor z=z_1 + i z_2 \in D(A_\C).
\end{align*}  
Then \emph{the complex resolvent set} and \emph{the real resolvent set} of the operator $A$ are given by
\begin{align*}
& \varrho(A) := \{\lambda\in \C \ | \ \ker (\lambda I - A_\C) = \{0\},\ (\lambda I - A_\C)^{-1}\in L(X_\C)\}, \\
& \varrho(A,\R) := \{\lambda\in \R \ | \ \ker (\lambda I - A) = \{0\}, \ (\lambda I - A)^{-1}\in L(X)\}
\end{align*}
respectively. Furthermore \emph{the complex spectrum} and {\em the real spectrum} of $A$ are defined by
\begin{align*}
\sigma(A):= \{\lambda\in\C \ | \ \lambda\not\in\varrho(A)\} \quad\text{and}\quad \sigma(A, \R):=\{\lambda\in\R \ | \ \lambda\not\in\varrho(A,\R)\},
\end{align*}
respectively, and we write
\begin{align*}
& \sigma_p(A) :=\{\lambda\in\C \ | \ \ker (\lambda I -A_\C) \neq \{0\}\}, \\
& \sigma_p(A,\R):= \{\lambda\in\R \ | \ \ker (\lambda I -A) \neq \{0\}\}
\end{align*}
for {\em the point spectrum} and {\em the real point spectrum}, respectively.
\begin{remark}\label{rem-pom}
The spectrum $\sigma(A)$ consists of the sequence (possibly finite) of real eigenvalues. Indeed, the operator $A$ has compact resolvents, and therefore $A_\C$ has also compact resolvents which implies that $$\sigma(A) = \sigma(A_\C,\C) = \sigma_p(A_\C,\C) = \{\lambda_i \ | \ i\ge 1\},$$ where $(\lambda_i)$ is finite or $|\lambda_i|\to +\infty$ when $n\to +\infty$. Furthermore, if $\lambda\in\C$ is an eigenvalue of $A_\C$, then it is also eigenvalue of the symmetric operator $\h A_\C$ and hence $\lambda$ is a real number. \hfill $\square$
\end{remark}
Let $Y\subset X$ be a linear subspace of $X$. \emph{The part of the operator $A$ in the space $Y$} is a linear operator $A_Y:Y\supset D(A_Y)\to Y$ given by
\begin{align}
D(A_Y) & :=\{x\in D(A) \ | \ Ax\in Y\}, \\
A_Y x & := Ax \ffor x\in D(A_Y).
\end{align}
We first prove the following theorem concerning spectral decomposition of $A$.
\begin{theorem}\label{th:10}
Assume that $(A1)$, $(A2)$ and $(A3)$ hold. If $\lambda = \lambda_k$ for some $k\ge 1$ is an eigenvalue of the operator $A$ and $X_0 := \ker (\lambda I - A)$, then $X = X_+\oplus X_-\oplus X_0$ for closed subspaces $X_+$, $X_-$ of $X$ and the following holds.
\begin{enumerate}
\item[(i)] One has $X_-\subset D(A)$, $A(X_-)\subset X_-$, $A(X_+\cap D(A)) \subset X_+$, $X_-$ is a finite dimensional space such that $X_-= \{0\}$ provided $k=1$ and
\begin{equation}\label{dir-sum}
 X_-=\bigoplus_{i=1}^{k-1} \ker(\lambda_i I - A) \text{ if } k\ge 2.
\end{equation}
Hence $\dim X_- = 0$ if $k=1$ and $\dim X_- = \sum_{i=1}^{k-1} \dim\ker(\lambda_i I - A)$ if $k\ge 2$.
\item[(ii)] If $A_+:X_+\supset D(A_+) \to X_+$ and $A_-:X_-\supset D(A_-) \to X_-$ are parts of $A$ in $X_+$ and $X_-$, respectively, then $\sigma(A_+) = \{\lambda_i \ | \ i\ge k+1 \}$ and
\begin{equation*}
 \sigma(A_-) = \begin{cases} \emptyset & \text{ if } k=1, \\ \{\lambda_i \ | \ i=1,\ldots, k-1 \}  & \text{ if } k\ge 2. \end{cases}
\end{equation*}
\item[(iii)] We have $\<i (u_l),i(u_m)\>_H = 0$ for $u_l\in X_l$ and $u_m\in X_m$ where $l,m\in\{0,-,+\}$, $l\neq m$.
\end{enumerate}
\end{theorem}
In the proof we use the following lemmata.
\begin{lemma}\label{lem-wart} 
Let $D:W\to W$ be a linear operator on a real finite dimensional space $V$ such that $W=W_1\oplus W_2\oplus\ldots\oplus W_l$ ($l\ge 1$) and $D x = \mu_i x$ for $x\in W_i$, where $\mu_i\in\R$ ($1\le i \le l$). Then
\begin{enumerate}
\item[(a)] $\sigma(D, \R) = \sigma(D) = \{\mu_i \ | \ 1\le i \le l\}$, \\[-10pt]
\item[(b)] For any $1\le i \le l$ we have $N_{\mu_i}(D) = \ker(\mu_i I - D)$.
\end{enumerate}
\end{lemma}
\noindent\textbf{Proof.} $(a)$ We prove that $\sigma(D) \subset \{\mu_i \ | \ 1\le i \le l\}$. The opposite inclusion is clear. Take $\mu\in\C$ such that
$\mu z = D_\C z$ where $z := x + iy\in W_\C $, $z\neq 0$. Then we have $W_\C = W_1\times W_1\oplus W_2\times W_2\oplus\ldots\oplus W_l\times W_l$ and $D_\C z = \mu_i z$ for $z\in W_i\times W_i$ ($1\le i \le l$). Hence $z=z_1+z_2+\ldots+ z_l$ where $z_i\in W_i\times W_i$ ($1\le i\le l$) and therefore
$\mu z = D_\C z =  \mu_1 z_1 + \mu_2 z_2 + \ldots + \mu_l z_l$. Since $z\neq 0$, there exists $1\le i \le l$ such that $z_i\neq 0$ and therefore $\mu = \mu_i$, which gives desired inclusion. \\[5pt]
$(b)$ It is enough to show that $N_{\mu_i}(D) \subset \ker(\mu_i I - D)$. If we take $x\in N_{\mu_i}(D)\setminus \{0\}$, then there exists $i_0 \ge 1$ such that $(\mu_i I - D)^{i_0} x = 0$ and $x_i\in W_i$ ($1\le i\le l$) such that $x = x_1 + x_2 + \ldots + x_l$. Therefore
\begin{align*}
0 & = (\mu_i I - D)^{i_0} x = (\mu_i I - D)^{i_0} x_1 + (\mu_i I - D)^{i_0} x_2 + \ldots + (\mu_i I - D)^{i_0} x_l \\
& = (\mu_i - \mu_1)^{i_0} x_1 + (\mu_i - \mu_2)^{i_0} x_2 + \ldots + (\mu_i - \mu_l)^{i_0} x_l.
\end{align*} 
Since $x \neq 0$, one of $x_1, x_2, \ldots, x_n$ has to be also nonzero. If we assume that $x_j \neq 0$ for some $1\le j \le l$, then $(\mu_i - \mu_j)^{i_0} x_j = 0$ and therefore $\mu_i = \mu_j$. This yields $x\in \ker(\mu_i I - B)$, which gives desired inclusion.  \hfill $\square$ \\

\begin{lemma}{\em (see \cite{MR610244})}\label{lem-rozkl-spect-com}
Let $A:X\supset D(A)\to X$ be a linear operator on a Banach space $X$ and assume that $X= X_-\oplus X_0\oplus X_+$ for closed subspaces $X_0$, $X_-$, $X_-$ such that
$$X_0,X_-\subset D(A), \ A(X_0)\subset X_0, \ A(X_-)\subset X_- \ \text{ and } \ A(D(A)\cap X_+)\subset X_+.$$
Let the operator $A_i : X_i\supset D(A_i) \to X_i$ be a part of the operator $A$ in the space $X_i$ for any $i=0,-,+$. Then the following assertions hold.
\begin{enumerate}
\item[(a)] For any $i=0,-,+$, we have $\varrho(A,\R)\subset \varrho(A_i,\R)$ and furthermore, if $\rho\in\varrho(A,\R)$ then
\begin{equation}\label{zw-re}
    (\rho I - A_i)^{-1} x  = (\rho I - A)^{-1}x \ffor x\in X_i.
\end{equation}
\item[(b)] If $A$ has compact resolvents, then for any $i=0,-,+$ the operator $A_i$ has also compact resolvents.
\item[(c)] If $-A$ is a generator of a $C_0$ semigroup $\{S_A(t)\}_{t\ge 0}$, then $$S_A(t) X_i \subset X_i \ffor t\ge 0 \ \text{ and } \ i=0,-,+.$$
\end{enumerate}
\end{lemma}

\begin{lemma}\label{part-op}
Under assumption $(A1)$, $(A2)$ and $(A3)$ the following assertions hold.
\begin{enumerate}
\item[(a)] For any $l\ge 1$ the following equality holds
\begin{equation}\label{eq-ker}
   \ker(\lambda_l I - A) = N_{\lambda_l}(A).
\end{equation}
\item[(b)] If $Y\subset X$ is a subspace of $X$, then $\sigma_p(A_Y) = \sigma_p(A_Y, \R)$.
\end{enumerate}
\end{lemma}

\noindent\textbf{Proof.} $(a)$ Since the operator $A$ has compact resolvents, there is $i_0 \ge 1$ such that $N_{\lambda_l}(A) = \ker(\lambda_l I - A)^{i_0}$. If we take $x\in\ker (\lambda_l I - A)^{i_0}$ then, from assumption $(A3)$, we have $(\lambda_l I - \h A)^{i_0} i(x) = 0$ and therefore $(\lambda_l I - \h A) i(x) = 0$, because the operator $\h A$ is symmetric. Consequently, $x \in \ker(\lambda_l I - A)$, which proves \eqref{eq-ker}. \\[5pt]
$(b)$ If $\lambda\in \sigma_p(A_Y)$ then $\lambda$ is an eigenvalue of the operator $(A_Y)_\C$, and hence is also an eigenvalue of the operator $A_\C$. Hence $\lambda$ is a real number as Remark \ref{rem-pom} says. Since $\sigma_p(A_Y)\cap \R = \sigma_p(A_Y,\R)$, it follows that $\lambda\in\sigma_p(A_Y,\R)$ and hence $\sigma_p(A_Y) \subset \sigma_p(A_Y, \R)$. The opposite inclusion is immediate. \hfill $\square$ \\

\begin{theorem}\label{zwa-rez}
Let $A:X\supset D(A)\to X$ be a linear operator with compact resolvents on a real Banach space $X$ and let $(\lambda_i)_{i\ge 1}$ be a sequence of real eigenvalues of the operator $A$. Then for any $k\ge 1$, there is a direct sum decomposition $X=X_1\oplus X_2$ such that $X_1$, $X_2$ are closed, $$X_1 = \bigoplus_{l=1}^k N_{\lambda_l}(A) \quad\text{ and }\quad  X_2 = \bigcap_{l=1}^k R_{\lambda_l}(A)$$
and the following assertions hold:
\begin{enumerate}
\item[(a)] $X_1\subset D(A)$, $A(X_1)\subset X_1$ and $A(X_2\cap D(A)) \subset X_2$,\\[-10pt]
\item[(b)] if $A_1$ and $A_2$ are parts of the operator $A$ in $X_1$ and $X_2$, respectively, then
$$\sigma(A_1, \R) = \{\lambda_1, \lambda_2, \ldots, \lambda_k\}\quad\text{and}\quad\sigma(A_2,\R) = \{\lambda_i \ | \ i\ge k + 1\}.$$
\end{enumerate}
\end{theorem}
\noindent\textbf{Proof.} Let $\rho\in\varrho(A, \R)$. Then $\lambda_i = \rho - \mu_i^{-1}$ for $i\ge 1$, where $\sigma_p((\rho I - A)^{-1}, \R) = \{\mu_i \ | \ i\ge 1\}$. Then (see \cite{MR2020421}) there is a direct sum decomposition $X=X_1\oplus X_2$ on closed subspaces such that
$$X_1 = \bigoplus_{l=1}^k N_{\mu_l}((\rho I - A)^{-1}) \quad\text{ and }\quad  X_2 = \bigcap_{l=1}^k R_{\mu_l}((\rho I - A)^{-1}).$$ Furthermore, $(\rho I - A)^{-1}(X_1)\subset X_1$, $(\rho I - A)^{-1}(X_2) \subset X_2$ and
\begin{equation}\label{row-spec2}
\sigma_p((\rho I - A)^{-1}_{|X_1}, \R) = \{\mu_1, \mu_2, \ldots, \mu_k\}, \quad \sigma_p((\rho I - A)^{-1}_{|X_2}, \R) = \{\mu_i \ | \ i\ge k + 1\}.
\end{equation}
On the other hand, one can check that $$N_{\lambda_i}(A) = N_{\mu_i}((\rho I - A)^{-1}) \quad\text{ and }\quad R_{\lambda_i}(A) = R_{\mu_i}((\rho I - A)^{-1}) \ffor i\ge 1,$$ which implies that $X=X_1\oplus X_2$ where $$X_1 = \bigoplus_{l=1}^k N_{\lambda_l}(A) \quad\text{ and }\quad  X_2 = \bigcap_{l=1}^k R_{\lambda_l}(A).$$
It is not difficult to verify that $$X_1\subset D(A), \ A(X_1)\subset X_1 \ \text{ and } \ A(D(A)\cap X_2)\subset X_2$$ and hence, by Lemma \ref{lem-rozkl-spect-com} $(a)$, we obtain
\begin{equation}\label{row-rez}
    \rho\in\varrho(A_i,\R) \quad\text{ and } \quad (\rho I - A_i)^{-1} = (\rho I - A)^{-1}_{|X_i} \ffor i=1,2.
\end{equation}
Further, from the point $(b)$ of the same lemma, it follows that the operators $A_1$ and $A_2$ have compact resolvents and therefore
\begin{equation*}
\sigma(A_i,\R) = \sigma_p(A_i,\R) = \{\rho - \mu^{-1} \ | \ \mu\in\sigma_p((\rho I - A_i)^{-1},\R)\} \ffor i=1,2,
\end{equation*}
which together with \eqref{row-spec2} and \eqref{row-rez} yields
\begin{equation*}
\begin{aligned}
    \sigma(A_1,\R) & = \{\rho - \mu_i^{-1} \ | \ 1\le i \le k\} = \{\lambda_i \ | \ 1\le i \le k\} \ \ \text{ and } \\
    \sigma(A_2,\R) & = \{\rho - \mu_i^{-1} \ | \ i\ge k+1\} = \{\lambda_i \ | \ i\ge k+1\},
\end{aligned}
\end{equation*}
and the proof is completed. \hfill $\square$ \\

\noindent\textbf{Proof of Theorem \ref{th:10}.} By Theorem \ref{zwa-rez} we obtain a direct sum decomposition of the space $X=X_-\oplus N_{\lambda_k}(A)\oplus X_+$, where $X_-= \{0\}$ if $k=1$, $$X_-=\bigoplus_{i=1}^{k-1} \ker(\lambda_i I - A) \text{ if } k\ge 2 \quad\text{ and } \quad  X_+ = \bigcap_{i=1}^k R_{\lambda_i}(A).$$ Furthermore we have the inclusions $X_-\subset D(A)$, $A(X_-)\subset X_-$, $A(X_+\cap D(A)) \subset X_+$ and
\begin{equation}\label{spek}
\sigma(A_+,\R) = \{\lambda_i \ | \ i\ge k+1 \}.
\end{equation}
From Lemma \ref{part-op} $(a)$ we infer that
\begin{equation}\label{rzb}
   \ker(\lambda_l I - A) = N_{\lambda_l}(A) \ffor l\ge 1.
\end{equation}
Hence $X_-$ is finite dimensional and
$\dim X_- = 0$ if $\lambda = \lambda_1$ and
$$\dim X_- =\sum_{i=1}^{k-1} \dim N_{\lambda_i}(A)  = \sum_{i=1}^{k-1} \dim \ker(\lambda_i I - A),$$ if $\lambda = \lambda_k$ for some $k\ge 2$. In this way we proved point $(i)$. \\[5pt]
In order to prove point $(ii)$, observe that from Lemma \ref{lem-rozkl-spect-com} $(b)$ it follows that the operator $A_+$ has compact resolvents. Hence the operator $(A_+)_\C$ also has compact resolvents and therefore $\sigma(A_+) = \sigma_p (A_+)$. By Lemma \ref{part-op} $(b)$ we infer that $\sigma_p(A_+) = \sigma_p(A_+,\R)$. From Lemma \ref{lem-rozkl-spect-com} $(b)$ it follows that the operators $A_+$ and $A_-$ have compact resolvents and hence $\sigma_p(A_+,\R) = \sigma(A_+,\R)$. This together with \eqref{spek} gives $\sigma(A_+) = \{\lambda_i \ | \ i\ge k+1 \}$.

If $k = 1$ then $X_- = \{0\}$ and hence $\sigma(A_-) = \emptyset$. If we suppose that $k\ge 2$, then combining the inclusion $A(X_-)\subset X_-$, \eqref{rzb} and Lemma \ref{lem-wart} we deduce that $\sigma(A_-) = \{\lambda_i \ | \ i=1,\ldots, k-1\}$ and the proof of point $(ii)$ is completed. \\[5pt]
We proceed to point $(iii)$. Take $1\le l \le k$ and $x\in N_{\lambda_l}(A)$, $y\in X_+$. Then $y\in R_{\lambda_l}(A)$ and furthermore, by \eqref{rzb}, we have $i(x)\in\ker(\lambda_l I - \h A)$ and $i(y)\in\im(\lambda_l I - \h A)$. Since the operator $\h A$ is symmetric, we have $\<i(x), i(y)\>_H = 0$. Accordingly, for any $1\le l\le k$ the spaces  $i(N_{\lambda_l}(A))$ and $i(X_+)$ are orthogonal, which implies that the spaces $i(X_+)$, $i(X_-)$ and $i(X_0)$, $i(X_+)$ are mutually orthogonal. Now we take $x\in N_{\lambda_k}(A)$ and $y\in N_{\lambda_l}(A)$, where $1\le l \le k-1$. In view of \eqref{rzb} we infer that $i(x)\in\ker(\lambda_k I - \h A)$ and $i(y)\in\ker(\lambda_l I - \h A)$, which along with the fact that $\h A$ is symmetric gives $\<i(x), i(y)\>_H = 0$. Hence the spaces $i(X_-)$ and $i(X_0)$ are also orthogonal and the proof of point $(iii)$ is completed.  \hfill $\square$ \\

\begin{lemma}{\em (see \cite[Theorem 1.5.3]{MR610244})}\label{lem-rozkl-spect}
Let $\lambda = \lambda_k$ for some $k\ge 1$ be an eigenvalue of the sectorial operator $A$ and let $X= X_1\oplus X_0\oplus X_2$ be a direct sum decomposition such that $X_0 = \lambda I - A$ and
$$X_0,X_-\subset D(A), \ A(X_0)\subset X_0, \ A(X_-)\subset X_- \ \text{ and } \ A(D(A)\cap X_+)\subset X_+.$$
Assume that $A_1:X_1\supset D(A_1) \to X_1$ and $A_2:X_2\supset D(A_2) \to X_2$ are parts of the operator $A$ in $X_1$ and $X_2$, respectively. If
$$\re z < 0 \ \text{ for } \ z\in\sigma(\lambda_k I - A_1)\quad \text{and}\quad\re z > 0 \ \text{ for } \ z\in\sigma(\lambda_k I - A_2),$$
then there are positive constants $c_\alpha$ and $C_\alpha$ such that
\begin{alignat}{2}\label{r1a}
\|A_\delta^\alpha S_A(t)x\| & \le C_\alpha t^{-\alpha} \, e^{- (\lambda + c_\alpha)t} \|x\|   && \ffor t > 0, \ x\in X_+, \\ \label{r3a}
\|S_A(t)x\| & \le C_\alpha e^{-(\lambda + c_\alpha) t}\|x\| && \ffor t\ge 0, \ x\in X_+ \\ \label{r2a}
\|S_A(t)x\| & \le C_\alpha e^{- (\lambda - c_\alpha)t} \|x\|   && \ffor t \le 0, \ x\in X_-,
\end{alignat}
where $S_A(t)x:= \exp(-t A_-)x$ for $t\in\R$ and $x\in X_-$ is the natural extension of the semigroup $\{S_A(t)\}_{t\ge 0}$ on the space $X_-$.
\end{lemma}

\noindent\textbf{Proof of Theorem \ref{th:10ab}.} By Theorem \ref{th:10} we obtain a direct sum decomposition $X=X_0\oplus X_-\oplus X_+$ such that assertions $(i)$ and $(iii)$ are satisfied. Furthermore
$$X_0,X_-\subset D(A), \ A(X_0)\subset X_0, \ A(X_-)\subset X_-, \ A(D(A)\cap X_+)\subset X_+$$
and from point $(ii)$ it follows that $$\re z < 0 \ \text{ for } \ z\in\sigma(\lambda_k I - A_+)\quad \text{and}\quad\re z > 0 \ \text{ for } \ z\in\sigma(\lambda_k I - A_-),$$ where $A_+:X_+\supset D(A_+) \to X_+$ and $A_-:X_-\supset D(A_-) \to X_-$ are parts of the operator $A$ in $X_+$ and $X_-$, respectively. Hence Lemma \ref{lem-rozkl-spect} implies point $(ii)$ while Lemma \ref{lem-rozkl-spect-com} $(c)$ leads to inclusions $$S_A(t)X_i \subset X_i \ffor t\ge 0, \ i\in\{0,-,+\}$$
and the proof is completed. \hfill $\square$ \\

\begin{theorem}{\em (see \cite[Theorem 16.7.2]{MR0089373})}\label{th-sem-spec}
Let $\{S_A(t)\}_{t\ge 0}$ be a $C_0$ semigroup on a Banach space $X$, which is generated by $- A$. Then \begin{equation*}
e^{-t\sigma_p(A)}\subset \sigma_p(S_A(t))\subset e^{-t\sigma_p(A)} \cup \{0\} \ffor t > 0.
\end{equation*}
Furthermore, if $\lambda\in\C$ then
\begin{equation}\label{th-sem-1}
\ker(e^{- \lambda t} I - S_A(t)) = \overline{\mathrm{span}}\left(\bigcup_{k\in\mathbb{Z}}\ker(\lambda_{k,t} I - A)\right),
\end{equation}
where $\lambda_{k,t} = \lambda + (2k\pi /t )i$ \ for \ $k\in\mathbb{Z}$.
\end{theorem}

\section{Cauchy problems}

Consider differential equations of the form
\begin{equation}\label{row-a-fs}
\dot u(t) = - A u(t) + \lambda u(t) + F (s,t,u(t)), \qquad  t > 0,
\end{equation}
where $s\in[0,1]$ is a parameter, $\lambda$ is a real number, $A: X\supset D(A)\to X$ is a sectorial operator with compact resolvents on a Banach space $X$ and $F:[0,+\infty)\times X^\alpha\to X$ is a continuous map. Here $X^\alpha$ for $\alpha\in(0,1)$, is a fractional power space given by $X^\alpha:=D((A + \delta I)^\alpha)$, where $\delta > 0$ is such that $A + \delta I$ is a positively defined operator. Furthermore we assume that: \\[5pt]
\noindent\makebox[9mm][l]{$(F1)$}\parbox[t]{118mm}{for any $x\in X^\alpha$ there exists a neighborhood $x\in V\subset X^\alpha$ and $L > 0$ such that for $s\in[0,1]$, $x_1,x_2\in V$ and $t\in [0,+\infty)$ one has $$\|F(s,t,x_1) - F(s,t,x_2)\|\le L \|x_1 - x_2\|_\alpha;$$}\\
\noindent\makebox[9mm][l]{$(F2)$}\parbox[t]{118mm}{for any $s\in[0,1]$, $t\in [0,+\infty)$, $x\in X^\alpha$ one has 
$$\|F(s,t,x)\| \le c(t)(1 + \|x\|_\alpha),$$ where $c:[0,+\infty)  \to [0,+\infty)$ is a continuous function.}\\

In this section we intend to recall the facts concerning existence and uniqueness of mild solutions for the equation \eqref{row-a-fs}. Then we provide theorems for the continuous dependence from parameter and initial data. Finally some compactness properties of mild solutions will be formulated.
\begin{definition}
Given an interval $J\subset\R$, we say that a continuous mapping $u:J \to X^\alpha$ is \emph{a mild solution} of equation \eqref{row-a-fs}, if
\begin{equation*}
u(t) = S_A(t - t')u(t') + \int_{t'}^t S_A(t - \tau)F(s,\tau,u(\tau)) \d \tau
\end{equation*}
for every $t,t' \in J$, $t' < t$.
\end{definition}
\begin{remark}\label{rem-klc}
Assume that the operator $A:X\to X$ is bounded and let $u:J \to X^\alpha$ be a mild solution of \eqref{row-a-fs}. \\
(a) Then it is known that $u$ is a $C^1$ map on $J$ and the equation \eqref{row-a-fs} is satisfied for any $t\in J$. \\
$(b)$ For any $t,t'\in J$ we have
\begin{equation}\label{row-po}
u(t) = S_A(t - t')u(t') + \int_{t'}^t S_A(t - \tau)F(s,\tau,u(\tau)) \d \tau
\end{equation}
where $S_A(t) := \exp(-tA)$ for $t\in \R$. Indeed, for $t,t' \in J$ such that $t' > t$ we have
\begin{equation*}
u(t') = S_A(t' - t)u(t) + \int_t^{t'} S_A(t' - \tau)F(s,\tau,u(\tau)) \d \tau.
\end{equation*}
Acting on this equation by $S_A(t - t')$ we derive
\begin{equation*}
S_A(t - t')u(t') = u(t) + \int_t^{t'} S_A(t - \tau)F(s,\tau,u(\tau)) \d \tau,
\end{equation*}
which implies \eqref{row-po}. \hfill $\square$
\end{remark}

\begin{theorem}{\em (see \cite[Theorem 3.3.3, Corollary 3.3.5]{MR610244})}\label{th-exist1}
For every $s\in[0,1]$ and $x\in X^\alpha$, the equation \eqref{row-a-fs} admits a unique mild solution $u(t;s, x):[0,+\infty)\to X^\alpha$ starting at $x$.
\end{theorem}

As we will use the topological invariants we need some continuity and compactness properties for the solutions. Here the key point is the assumption concerning the compactness of resolvents of the operatora $A$, which will be used to prove the following theorems.

\begin{theorem}\label{tw-con-21}
If $(x_n)$ in $X^\alpha$ and $(s_n)$ in $[0,1]$ are such that $x_n\to x_0$ in $X^\alpha$ and $s_n\to s_0$ when $n\to +\infty$, then
\begin{equation}\label{rrr111}
u(t;s_n, x_n) \to u(t;s_0, x_0) \aas n\to +\infty,
\end{equation}
for any $t\ge 0$, and furthermore this convergence is uniform for $t$ from bounded subsets of $[0,+\infty)$.
\end{theorem}

\begin{theorem}\label{tw-exi-con-comp}
If $t > 0$ and $\Omega\subset X^\alpha$ is a bounded set, then $$\{u(t;s,x) \ | \ s\in[0,1], \ x\in\Omega\}$$ is a relatively compact subset of $X^\alpha$.
\end{theorem}

Before we strat the proof we formulate some auxiliary lemmata.

\begin{lemma}{\em (see \cite[Lemma 1.2.9]{MR1778284})} \label{lem-volt-type-ineq}
Let $\alpha\in [0,1)$, $a\ge 0$, $b > 0$ and let $\phi\colon [t_0,T)\to [0,+\infty)$ be a continuous function such that
\begin{equation*}
\phi(t)\le a + b\int_{t_0}^t \frac{1}{(t - \tau)^\alpha} \phi(\tau) \d \tau \ffor  t\in (t_0, T).
\end{equation*}
Then $$\sup_{t\in [t_0, T)}\phi(t)\le a K(\alpha,b,T),$$ where $K(\alpha,b,T)$ is a constant dependent from $\alpha$, $b$ and $T$.
\end{lemma}

\begin{lemma}\label{rem-bound-f}
Let $(s_n)$ be a sequence in $[0,1]$ and let $(x_n)$ be a bounded sequence in $X^\alpha$. Assume that, for any $n\ge 1$, the map $v_n:[0,t_0]\to X$ is given by $$v_n(\tau) := F(s_n,\tau,u(\tau;s_n, x_n)) \ffor \tau\in [0,t_0].$$ Then the following assertions hold.
\begin{enumerate}
\item[(a)] The set $\{u(t ;s_n, x_n) \ | \ t\in [0,t_0], \  n\ge 1\}$ is bounded in $X^\alpha$.
\item[(b)] For any $t\in [0,t_0]$, the set
$$\left\{\int_0^t A_\delta^\alpha S_A(t - \tau)v_n(\tau)\d \tau \ \Big| \ n\ge 1\right\}$$
is bounded in $X$.
\item[(c)] For any $\varepsilon > 0$, there is $h_0 > 0$ such that, if $t,t'\in [0,t_0]$ and $0 < t' - t < h_0$, then
$$\int_t^{t'} \|A_\delta^\alpha S_A(t' - \tau)v_n(\tau)\| \d \tau  \le \varepsilon \ffor n\ge 1.$$
\end{enumerate}
\end{lemma}
\noindent\textbf{Proof.} Let $c\in\R$, $M$ ,$M_\alpha$ be a constants such that $$\|S_A(t)\| \le M \ \ \text{for} \ \ t\in [0,t_0] \ \text{ and } \ \|A_\delta^\alpha S_A(t)\| \le M_\alpha e^{c t} t^{-\alpha} \ \ \text{for} \ \ t > 0.$$
Since $(x_n)$ is bounded, there is $R > 0$ such that $\|x_n\|_\alpha \le R$ for $n\ge 1$. Then, assumption $(F2)$ implies that for any $n\ge 1$ and $t\in[0,t_0]$, we have
\begin{align*}
\|u(t;s_n, x_n)\|_\alpha & \le  \|S_A(t) A_\delta^\alpha x_n\|
+ \int_0^t \|A_\delta^\alpha S_A(t - \tau)F(s_n,\tau,u(\tau;s_n, x_n))\|\d \tau \\
& \le M \|x_n\|_\alpha + \int_0^t \frac{M_\alpha e^{ c (t - \tau)} }{(t - \tau)^\alpha}\|F(s_n,\tau,u(\tau;s_n, x_n))\|\d \tau \\
& \le M R + \int_0^{t} \frac{M_\alpha e^{| c| t_0}}{(t - \tau)^\alpha} \, c(\tau)(1 + \|u(\tau;s_n, x_n)\|_\alpha) \d \tau \\
& \le M R + \frac{K M_\alpha e^{| c| t_0}}{1-\alpha} t_{0}^{1-\alpha}
+ \int_0^{t} \frac{K M_\alpha e^{| c| t_0}}{(t - \tau)^\alpha}\,\|u(\tau;s_n, x_n)\|_\alpha \d \tau,
\end{align*}
where $K := \sup_{\tau\in [0,t_0]}c(\tau)$. Hence, by Lemma \ref{lem-volt-type-ineq}, there is $C > 0$ such that $\|u(t;x_n,s_n)\|_\alpha \le C$ for $t\in [0,t_0]$ and $n\ge 1$, which proves $(a)$. Furthermore, note that by assumption $(F2)$ we have
\begin{align*}
\|v_n(t)\| = \|F(s_n,\tau,u(t;s_n, x_n))\|
\le c(t)(1 + \|u(t;s_n, x_n)\|_\alpha) \le K(1 + C)
\end{align*}
for $t\in [0,t_0]$ and $n\ge 1$. Then, for any $n\ge 1$, we infer that
\begin{align*}
\left\| \int_0^t A_\delta^\alpha S_A(t - \tau)v_n(\tau)\d \tau \right\|
& \le \int_0^t \|A_\delta^\alpha S_A(t - \tau)v_n(\tau)\|\d \tau \\
& \le \int_0^t \frac{M_\alpha e^{ c (t-\tau)}}{(t - \tau)^\alpha}\|v_n(\tau)\|\d \tau \\
& \le \int_0^t K(1 + C) \frac{M_\alpha e^{| c| t_0}}{(t - \tau)^\alpha} \d \tau
\le K(1 + C) \frac{M_\alpha e^{| c| t_0} }{1 - \alpha}{t_0}^{1 - \alpha},
\end{align*}
which gives $(b)$. As for $(c)$, let $t',t\in [0,t_0]$ be such that $t' > t$. Then, for any $n\ge 1$, we obtain
\begin{align*}
\left\| \int_t ^{t'} A_\delta^\alpha S_A(t' - \tau)v_n(\tau)\d \tau \right\|
& \le \int_t ^{t'} \frac{M_\alpha e^{ c (t' - \tau)}}{(t' - \tau)^\alpha}\|v_n(\tau)\|\d \tau \\
& \le \int_t ^{t'} K(1 + C)\frac{M_\alpha e^{| c| t_0}}{(t' - \tau)^\alpha}
= K(1 + C) \frac{M_\alpha e^{| c| t_0}}{1 - \alpha}(t' - t)^{1 - \alpha}.
\end{align*}
Take $h_0 := \left(\frac{\varepsilon (1 - \alpha)}{K(1 + C)M_\alpha e^{ c t_0}}\right)^{1/(1 - \alpha)}$. Then we see that, for any $t,t'\in [0,t_0]$ such that $0 < t' - t < h_0$, we have
$$\left\| \int_t^{t'} A_\delta^\alpha S_A(t' - \tau)v_n(\tau)\d \tau \right\| \le \varepsilon \ffor n\ge 1,$$
and the proof of point $(c)$ is completed. \hfill $\square$ \\

\noindent\textbf{Proof of Theorem \ref{tw-exi-con-comp}.} Let $t > 0$ and let $\Omega\subset X^\alpha$ be a bounded set. To prove the set $\Phi_t([0,1]\times \Omega )$ is relatively compact in $X^\alpha$ it is enough to prove that the set $A_\delta^\alpha\Phi_t([0,1]\times  \Omega)$ is relatively compact in $X$. To this end take sequences $(s_n)$ in $[0,1]$ and $(x_n)$ in $\Omega$ and let $(v_n)$ be given, for any $n\ge 1$, by $$v_n(\tau) := F(s_n,\tau,u(\tau ;s_n,x_n)) \ffor \tau\in [0,t].$$ If $\varepsilon > 0$ is arbitrary, then by Lemma \ref{rem-bound-f} $(c)$, there is $t_0\in (0,t)$ such that
\begin{equation*}
\left\|\int_{t_0}^t A_\delta^\alpha S_A(t - \tau)v_n(\tau) \d \tau \right\| \le \varepsilon \ffor n\ge 1.
\end{equation*}
Furthermore, by Lemma \ref{rem-bound-f} $(b)$, we infer that the set
\begin{equation*}
D_{t_0} := \left\{\int_0^{t_0} A_\delta^\alpha S_A(t_0 - \tau)v_n(\tau) \d \tau \ \Big| \ n \ge 1\right\}
\end{equation*}
is bounded. On the other hand, for any $n\ge 1$, we have
\begin{align*}
A_\delta^\alpha u(t;s_n, x_n) & = S_A(t)A_\delta^\alpha x_n
+ \int_{t_0}^t A_\delta^\alpha S_A(t - \tau)v_n(\tau)\d \tau \\
& \qquad + S_A(t - t_0)\left(\int_0^{t_0} A_\delta^\alpha S_A(t_0 - \tau)v_n(\tau) \d \tau \right),
\end{align*}
which implies that
\begin{align*}
V:=\{A_\delta^\alpha u(t ;s_n, x_n) \ | \ n\ge 1 \}
\subset S_A(t)\{A_\delta^\alpha x_n \ | \ n\ge 1\} + S_A(t - t_0)D_{t_0} \\
\quad + \left\{\int_{t_0}^t A_\delta^\alpha S_A(t_0 - \tau)v_n(\tau) \d \tau \ \Big| \ n \ge 1\right\}
\subset W + B(0,\varepsilon),
\end{align*}
where
$$W := S_A(t)\{A_\delta^\alpha x_n \ | \ n\ge 1\} + S_A(t - t_0)D_{t_0}.$$
Since the semigroup $\{S_A(t)\}_{t\ge 0}$ is compact and the sets $\{A_\delta^\alpha x_n \ | \ n\ge 1\}$, $D_{t_0}$ are bounded, w infer that the set $W$ is relatively compact in $X$. Since $\varepsilon > 0$ may be arbitrary small, we deduce that the set $V$ is also relatively compact in $X$, which completes the proof. \hfill $\square$\\

\begin{lemma}\label{lem-eq-con}
The family $\{u_n(\,\cdot\, ;s_n, x_n) \ | \ n\ge 1\}$ is equicontinuous for $t\in[0,+\infty)$.
\end{lemma}
\noindent\textbf{Proof.} For any $n\ge 1$ write $u_n := u_n(\,\cdot\, ;s_n, x_n)$. By the integral formula, for $t\in[0,+\infty)$, $h \ge 0$ and $n\ge 1$, we have
\begin{equation}
\begin{aligned}\label{est-theo-cont}
\|u_n(t+h) - u_n(t)\|_\alpha & \le \|S_A(h)u_n(t) - u_n(t)\|_\alpha  \\
& \quad + \int_t^{t+h} \|A_\delta^\alpha S_A(t + h - \tau)F(s_n, \tau,u_n(\tau)) \| \d \tau.
\end{aligned}
\end{equation}
Theorem \ref{tw-exi-con-comp} says that, for $\tau\in [0,t]$, the set $\{u_n(\tau) \ | \ n\ge 1\}$ is relatively compact in $X^\alpha$ and therefore there is $h_0 > 0$ such that
\begin{equation}\label{R1}
\|S_A(h)u_n(t) - u_n(t)\|_\alpha \le \varepsilon /2 \ffor  0< h < h_0, \quad n \ge 1.
\end{equation}
By Lemma \ref{rem-bound-f} $(c)$, we find that there is $h_1 > 0$ such that
\begin{equation}\label{R2}
\int_t^{t+h} \|A_\delta^\alpha S_A(t + h - \tau)F(s_n, \tau,u_n(\tau)) \| \d \tau \le\varepsilon /2
\end{equation}
for $0 < h <h_1$ and $n\ge 1$.
Combining \eqref{est-theo-cont}, \eqref{R1} and \eqref{R2} we deduce that
 $$\|u_n(t + h) - u_n(t)\|_\alpha \le \varepsilon /2 + \varepsilon /2 = \varepsilon \ffor 0< h < h_1,\quad n\ge 1$$
which implies that the family $\{u_n\}_{n\ge 1}$ is equicontinuous from the right side on $[0,+\infty)$. It remain to prove that the family is equicontinuous from the left side $(0,+\infty)$. To this end take $t\in (0,+\infty)$ and let $\varepsilon > 0$ be arbitrary. If $h$ and $\theta$ are such that $0 < h < \theta < t$, then
\begin{align*}
\|u_n(t) - u_n(t - h)\|_\alpha & \le \|u_n(t) - S_A(\theta)u_n(t - \theta)\|_\alpha \\
& \quad + \|S_A(\theta)u_n(t - \theta) - S_A(\theta - h)u_n(t - \theta)\|_\alpha \\
& \quad + \|S_A(\theta - h)u_n(t - \theta) - u_n(t - h)\|_\alpha
\end{align*}
and hence, for any $n\ge 1$, we have
\begin{equation}
\begin{aligned}\label{R6}
\|u_n(t) - u_n(t - h)\|_\alpha & \le
\int_{t - \theta}^t \|A_\delta^\alpha S_A(t - \tau)F(s_n, \tau,u_n(\tau)) \| \d \tau \\
& \quad + \|S_A(\theta)u_n(t - \theta) - S_A(\theta - h)u_n(t - \theta)\|_\alpha \\
& \quad + \int_{t - \theta}^{t - h} \|A_\delta^\alpha S_A(t - h - \tau)F(s_n, \tau,u_n(\tau))\|\d \tau.
\end{aligned}
\end{equation}
By Lemma \ref{rem-bound-f} $(c)$, there is $h_0 \in (0,t)$ such that, for any $t_1,t_2\in [0,t]$ with $0 < t_1 - t_2 < h_0$, we have
\begin{equation}\label{R5}
\int_{t_2}^{t_1} \|A_\delta^\alpha S_A(t_1 - \tau)F(s_n, \tau,u_n(\tau)) \| \d \tau \le \varepsilon /3 \ffor  n\ge 1.
\end{equation}
Let $\theta\in (0,h_0)$ be fixed. By Theorem \ref{tw-exi-con-comp}, the set $\{u_n(t - \theta) \ | \ n\ge 1\}$ is relatively compact and hence we can choose $h_1$ such that $0 < h_1 < \theta$ and
\begin{equation}\label{row-equicont}
\|S_A(\theta)u_n(t - \theta) - S_A(\theta - h)u_n(t - \theta)\|_\alpha \le \varepsilon /3 \quad\mbox{ for }\quad  h\in (0, h_1), \ n\ge 1.
\end{equation}
Using \eqref{R5}, for $h\in (0,h_1)$, we obtain
\begin{align}\label{R3}
\int_{t - \theta}^t \|A_\delta^\alpha S_A(t - \tau)F(s_n, \tau,u_n(\tau))\|\d \tau & \le \varepsilon /3 \qquad\mbox{oraz} \\ \label{R4}
\int_{t - \theta}^{t - h} \|A_\delta^\alpha S_A(t + h - \tau)F(s_n, \tau,u_n(\tau))\|\d \tau & \le \varepsilon /3 \ffor  n\ge 1.
\end{align}
Therefore, combining \eqref{R6}, \eqref{row-equicont}, \eqref{R3} and \eqref{R4} we infer that, for $h\in (0,h_1)$
\begin{equation*}
\|u_n(t) - u_n(t - h)\|_\alpha  \le \varepsilon /3 + \varepsilon /3 +\varepsilon /3 = \varepsilon,
\end{equation*}
and consequently the family $\{u_n \ | \ n\ge 1\}$ is equicontinuous from the left side on $(0,+\infty)$ as desired and the proof is completed. \hfill $\square$ \\

\noindent\textbf{Proof of Theorem \ref{tw-con-21}.} Write $u_n := u(\,\cdot\,;s_n, x_n)$ for $n\ge 1$. In view of Lemma \ref{lem-eq-con} and Theorem \ref{tw-exi-con-comp} we infer that, for any $T>0$, the family $(u_n)$ in equicontinuous and has relatively compact orbits on $[0,T]$. Let $(u_{n_k})_{k\ge 1}$ be arbitrary subsequence of $(u_n)_{n\ge 1}$. By Ascoli--Arzela Theorem theorem there is a subsequence $(u_{n_{k_l}})_{l\ge 1}$ and a continuous map $v:[0,T]\to X^\alpha$ such that $u_{n_{k_l}}(t)\to v(t)$ in $X^\alpha$, uniformly for $t\in [0,T]$ as $l\to +\infty$. Hence, letting $l\to +\infty$ in the formula
\begin{equation*}
u_{n_{k_l}}(t) = S_A(t) x_{n_{k_l}} + \int_0^t S_A(t - \tau)F(s_{n_{k_l}},\tau,u_{n_{k_l}}(\tau))\d \tau
\end{equation*}
for any $t\in[0,T]$, we have
\begin{equation*}
v(t) = S_A(t) x_0 + \int_0^t S_A(t - \tau)F(s_0,\tau,v(\tau))\d \tau.
\end{equation*}
Hence, by the uniqueness of mild solutions (see Theorem \ref{th-exist1}) we infer that $v(t) = u(t;s_0,x_0)$ for $t\in[0,T]$. Therefore $u_{n_{k_l}}(t)\to u(t;s_0,x_0)$ in $X^\alpha$, uniformly for $t\in [0,T]$ as $l\to +\infty$. Since the sequence $(u_{n_k})_{k\ge 1}$ is arbitrary, it follows that $u_n(t)\to u(t;s_0,x_0)$ in $X^\alpha$, uniformly for $t\in [0,T]$ as $n\to +\infty$ and the proof is completed. \hfill $\square$

\section{Resonant averaging principle}

We consider differential equations of the form
\begin{equation}\label{równanie-gl-eps}
\dot u(t) = - A u(t) + \lambda u(t) + \varepsilon F (t,u(t)), \qquad  t > 0.
\end{equation}
where $\lambda$ is an eigenvalue of $A$, and $F\colon [0,+\infty)\times X^\alpha \to X$ is a continuous map. Assume that $A$ and $F$ satisfies assumptions $(A1)$, $(A2)$, $(A3)$, $(F1)$ and \\[5pt]
\noindent\makebox[9mm][l]{$(F3)$}\parbox[t][][t]{118mm}{there is $m > 0$ such that $\|F(t,x)\|\le m$ for $t\ge 0$, $x\in X^\alpha$,}\\[5pt]
\noindent\makebox[9mm][l]{$(F4)$}\parbox[t][][t]{118mm}{there is $T > 0$ such that $F(t,x) = F(t + T, x)$ for $t\ge 0$, $x\in X^\alpha$.} \\[5pt]
By Theorem \ref{th-exist1}, the above assumptions imply that, for any $x\in X^\alpha$ and $\ve \ge 0$, there is a mild solution $u(\,\cdot\,;\varepsilon, x):\R\to X^\alpha$ of \eqref{równanie-gl-eps} starting at $x$. Let $\Phi_T:[0,1]\times X^\alpha \to X^\alpha$ be {\em the translation along trajectories operator} associated with this equation, given by $$\Phi_T(\varepsilon, x) := u(T;\varepsilon, x) \ffor \ve\ge 0, \ x\in X^\alpha.$$ Then, Theorems \ref{tw-con-21} and \ref{tw-exi-con-comp} say that $\Phi_T$ is a completely continuous map.

Remark \ref{rem-pom} says that the spectrum $\sigma(A)$ of the operator $A$ consists of the sequence of eigenvalues $$\lambda_1 < \lambda_2 < \ldots < \lambda_i < \lambda_{i+1} < \ldots$$ which is finite or $\lambda_i \to +\infty$ when $i\to +\infty$. Consider the direct sum decomposition $X=X_0\oplus X_-\oplus X_+$ on closed subspaces obtained in Theorem \ref{th:10ab}. Then $X_0 := \ker (\lambda I - A)$ and
\begin{equation}\label{ddir-sum}
    X_- = \bigoplus_{i=1}^{k-1} \ker(\lambda_i I - A).
\end{equation}
In particular $X_-$ is a finite dimensional space such that
\begin{equation}\label{wymia}
    \dim X_- = 0 \text{ if } k=1, \ \ \dim X_- = \sum_{i=1}^{k-1} \dim\ker(\lambda_i I - A) \text{ if } k\ge 2.
\end{equation}
It is also known that
\begin{equation}\label{inkl-sem}
S_A(t)X_i \subset X_i \ffor t\ge 0, \ i\in\{0,-,+\}.
\end{equation}
and there are constants $c, K > 0$ such that
\begin{align}\label{ine11}
     & \|A_\delta^\alpha S_A(t)x\| \le K e^{- (\lambda + c) t} t^{-\alpha} \|x\| && \hspace{-20mm}\text{for} \ x\in X_+, \ t > 0,\\ \label{ine22}
     & \|e^{\lambda t}S_A(t)x\| \le K e^{- c t}\|x\| && \hspace{-20mm}\text{for} \ x\in X_+, \ t\ge 0, \\ \label{ine33}
     & \|e^{\lambda t}S_A(t)x\| \le K e^{c t}\|x\| && \hspace{-20mm}\text{for} \ x\in X_-, \ t\le 0. 
\end{align}
where $S_A(t)x:= \exp(-t A_-)x$ for $t\in\R$ and $x\in X_-$ is the natural extension of the semigroup $\{S_A(t)\}_{t\ge 0}$ on $X_-$. Furthermore the spaces $X_0$, $X_-$ and $X_+$ are mutually orthogonal, that is,
\begin{equation}\label{ort}
    \<i (u_l),i(u_m)\>_H = 0 \ \ \text{for} \ \ u_l\in X_l, \ u_m\in X_m \text{ where } l,m\in\{0,-,+\}, \ l\neq m
\end{equation}
Let $P, Q_\pm:X\to X$ be projections given for any $x\in X$ by
\begin{equation}\label{wzz1ab}
    P x = x_0 \ \text{ and } \ Q_\pm x = x_\pm
\end{equation} 
where $x = x_+ + x_0 + x_-$ for $x_i\in X_i$, $i\in \{0,-,+\}$. Write $Q:= Q_- + Q_+$. Since the inclusion $X^\alpha\subset X$ is continuous, one can decompose $X^\alpha$ on a direct sum of closed spaces $X^\alpha= X_0\oplus X^\alpha_-\oplus X^\alpha_+$, where $X^\alpha_+:=X^\alpha\cap X_+$ and $X^\alpha_-:= X^\alpha\cap X_-$. Therefore the projections $P$ and $Q_\pm$ can be also considered as continuous maps $P, Q_\pm:X^\alpha\to X^\alpha$ given for any $x\in X^\alpha$ by \eqref{wzz1ab}. Note that by \eqref{inkl-sem}, we have
\begin{align}\label{ds1}
S_A(t)Px = PS_A(t)x \ \text{ and } \ S_A(t)Q_\pm x = Q_\pm S_A(t) x \ffor t\ge 0, \ x\in X.
\end{align}
Furthermore, by Theorem \ref{th-sem-spec} and Remark \ref{rem-pom}, we infer that
\begin{equation}\label{row-sp}
    \mathrm{Ker}\, (A - \lambda_i I) = \mathrm{Ker}\, (I - e^{\lambda_i t}S_A(t)) \ffor i\ge 1, \ t > 0.
\end{equation}
In this section we prove the following {\em resonant averaging principle}, expressing, for the small $\ve$, the Leray-Schauder topological degree $\mathrm{deg_{LS}}$ of the translation operator $\Phi_T(\ve,\,\cdot\,)$ in terms of the Brouwer degree $\mathrm{deg_B}$ of an appropriate averaging of the right side of \eqref{równanie-gl-eps}.
\begin{theorem}\label{th-aver-ker-sec}
Let $\lambda = \lambda_k$ for some $k\ge 1$ and let $h\colon N_0 \to N_0$ where $N_0:= \ker (\lambda I - A)$, be a map given by
\begin{equation*}
h(x) := \int_0^T PF(\tau,x)\,d \tau \qquad\mathrm{for}\quad x\in N_0.
\end{equation*}
Assume that $U \subset N_0$ and $V \subset X^\alpha_+\oplus X^\alpha_-$, where $0\in V$, are open and bounded subsets. If $h(x)\neq 0$ for $x \in \partial_{N_0} U$, then there is $\varepsilon_0 \in (0,1)$ such that for any $\varepsilon\in (0,\varepsilon_0]$ and $x\in \partial(U \oplus V)$ we have $\Phi_T(\varepsilon, x)\neq x$ and
\begin{equation*}
\mathrm{deg_{LS}}(I - \Phi_T(\varepsilon,\,\cdot\,), U \oplus V) = (-1)^{d_k} \cdot \mathrm{deg_B}(g,U),
\end{equation*}
where $d_0 := 0$ and $d_l:=\sum_{i=1}^l \dim\ker(\lambda_i I - A)$ for $l\ge 1$.
\end{theorem}
Let $\degb$ denote the Brouwer topological degree. In the proof we use the following theorem and lemma.
\begin{theorem}{\em (see \cite[Lemma 13.1]{MR736839})}\label{th-kras}
Consider the following differential equation $$\dot u(t) = \lambda f(u(t)), \qquad t > 0$$ where $\lambda\in[0,1]$ is a parameter and $f:\R^n \to \R^n$ is a bounded and continuous map. Let $\Theta^\lambda_T:\R^n\to\R^n$ be the translation operator associated with this equation. If $U\subset \R^n$ is an open bounded set such that $f(x)\neq 0$ for $x\in\partial U$, there is $\lambda_0 > 0$ such that, for $\lambda\in(0,\lambda_0]$ we have $\Theta^\lambda_T(x) \neq x$ and
\begin{equation*}
    \degb(I - \Theta^\lambda_T, U) = \degb(-f, U).
\end{equation*}
\end{theorem}

\begin{lemma}\label{lem-deg1}
If $\lambda = \lambda_k$ for some $k\ge 1$, is an eigenvalue of $A$, then
\begin{enumerate}
\item[(a)] $e^{\lambda T}S_A(T)x \neq x$ for $x\in X^\alpha_-\oplus X^\alpha_+$, $x\neq 0$, \\[-9pt]
\item[(b)] for any open set $V\subset X^\alpha_-\oplus X^\alpha_+$ such that $0\in V$ we have $$\degl(I - e^{\lambda T}S_A(T)|_{X^\alpha_-\oplus X^\alpha_+}, V) = (-1)^{d_{k-1}},$$
    where $d_0 := 0$ and $d_l:=\sum_{i=1}^l \dim\ker(\lambda_i I - A)$ for $l\ge 1$.
\end{enumerate}
\end{lemma}
\noindent\textbf{Proof.} $(a)$ Assume that $H:[0,1]\times X^\alpha_-\oplus X^\alpha_+ \to X^\alpha_-\oplus X^\alpha_+$ is a map given by $$H(\mu,x):= \mu e^{\lambda T}S_A(T) x_+ + e^{\lambda T}S_A(T) x_- \ffor x\in X^\alpha_-\oplus X^\alpha_+,$$ where for $x\in X^\alpha_-\oplus X^\alpha_+$ the elements $x_\pm\in X^\alpha_\pm$ are such that $x= x_+ + x_-$.
To prove the point $(a)$ we show that $H(\mu,x)\neq x$ for $\mu\in[0,1]$ and $x\in X^\alpha_-\oplus X^\alpha_+$ such that $x \neq 0$. Suppose by contradiction that $$\mu e^{\lambda T} S_A(T) x_+ + e^{\lambda T}S_A(T) x_- = x$$ for some $\mu\in[0,1]$ and $x\in X^\alpha_-\oplus X^\alpha_+$ such that $x \neq 0$. In view of \eqref{inkl-sem}, it implies that $\mu e^{\lambda T} S_A(T) x_+ = x_+$ and $e^{\lambda T} S_A(T) x_- = x_-$. We show that $x_+ = 0$. If $\mu = 0$ the it is immediate. If $\mu\in(0,1]$ then $S_{A - \lambda I}(T) x_+ = (1/\mu) x_+$, which along with Theorem \ref{th-sem-spec} and Remark \ref{rem-pom} implies that $x_+\in \ker((\lambda - \ln (1/\mu)/T) I - A)$. Since $\ln (1/\mu) \ge 0$, we infer that there is $1\le i \le k$ such that $A x_+ = \lambda_i x_+$. In view of \eqref{ddir-sum} it follows that $x_+\in X_-\oplus X_0$ and finally that $x_+ = 0$, because $x_+\in X_+$. Combining \eqref{row-sp} and the equality $e^{\lambda T} S_A(T) x_- = x_-$, we deduce that $x_-\in \ker(\lambda_k I - A) = X_0$. Since $x_-\in X_-$, we have $x_- = 0$ and hence $x = x_+ + x_- = 0$. This is impossible because we assumed that $x\neq 0$.  \\[5pt]
$(b)$ Let $V\subset X^\alpha_-\oplus X^\alpha_+$ be an open set such that $0\in V$. Since $H$ is an admissible homotopy, by the homotopy invariance of topological degree
\begin{align*}
\degl(I - e^{\lambda T} S_A(T), V) & = \degl(I - H(1,\,\cdot\,), V) = \degl(I - H(0,\,\cdot\,), V) \\
& = \degl(I - e^{\lambda T} S_A(T)|_{X_-}, V\cap X_-).
\end{align*}
By \eqref{ddir-sum} and the inclusion $\ker(\lambda_i I - A)\subset \ker(e^{(\lambda - \lambda_i) T}I - e^{\lambda T}S_A(T))$, Lemma \ref{lem-wart} implies that $$\sigma(e^{\lambda T} S_A(T)_{|X_-}, \R) = \{e^{(\lambda - \lambda_i) T} \ | \ 1\le i \le k-1\}$$ and the algebraic multiplicity of the eigenvalue $e^{(\lambda - \lambda_i) T}$, where $1\le i \le k-1$ is equal to $\dim \ker(\lambda_i I - A)$. Therefore we find that
\begin{align*}
\degl(I - e^{\lambda T} S_A(T), V) = \degl(I - e^{\lambda T} S_A(T)_{|X_-}, V\cap X_-) = (-1)^{d_{k-1}}
\end{align*}
which completes the proof. \hfill $\square$ \\

\noindent\textbf{Proof of Theorem \ref{th-aver-ker-sec}.} Consider the following family of differential equations
\begin{equation}\label{A-eps-res4}
\dot u(t) = - A u(t) + \lambda u(t) + \varepsilon G(s,t,u(t)), \qquad  t > 0
\end{equation}
where $G:[0,1]\times[0,+\infty)\times X^\alpha \to X$ is a map given by
$$G(s,t,x):= s F(t,x) + (1 - s) \frac{1}{T}\int_0^T PF(\tau,P x)\,d \tau$$ for $s\in[0,1]$, $t\in[0,+\infty)$, $x\in X^\alpha$. It is not difficult to check that $G$ satisfies assumptions $(F1)$ and $(F2)$, and hence, by Theorem \ref{th-exist1}, for any $x\in X^\alpha$ there is a mild solution $u(\,\cdot\,;s,\varepsilon, x):\R\to X^\alpha$ of \eqref{A-eps-res4} starting at $x$. Let $\Psi_T:[0,1]\times[0,1]\times X^\alpha \to X^\alpha$ be a translation operator associated with this equation given by $$\Psi_T(s,\varepsilon, x) := u(T;s,\varepsilon, x) \ffor (\varepsilon,s)\in [0,1]\times [0,1], \ x\in X^\alpha.$$
By Theorems \ref{tw-con-21} and \ref{tw-exi-con-comp} we infer that $\Psi_T$ is completely continuous. We show that there is $\varepsilon_0 > 0$ such that for any $\varepsilon\in (0,\varepsilon_0]$ we have $\Psi_T(\mu,\varepsilon,x)\neq x$ for $s\in[0,1]$ and $x\in \partial(U \oplus V)$. Otherwise there are sequences $(\varepsilon_n)$ in $(0,1]$, $(s_n)$ in $[0,1]$ and $(x_n)$ in $\partial(U \oplus V)$ such that $\varepsilon_n \to 0$ as $n\to +\infty$ and
\begin{equation}\label{c9}
\Psi_T(s_n,\varepsilon_n, x_n) = x_n \ffor n\ge 1.
\end{equation}
Since the operator $\Psi_T$ is completely continuous, by \eqref{c9}, we deduce that the sequence $(x_n)$ is relatively compact in $X^\alpha$. Therefore, without loss of generality we can assume that $s_n\to s_0$ and $x_n\to x_0$ as $n\to+\infty$, where $s_0\in[0,1]$ and $x_0\in \partial(U \oplus V)$. Letting $n\to +\infty$ in \eqref{c9} we have
\begin{equation*}
e^{\lambda T} S_A(T)x_0 = \Psi_T(0,0,x_0) = x_0,
\end{equation*}
which together with \eqref{row-sp} implies that
\begin{equation}\label{c11}
x_0 \in \ker(\lambda I - A) = N_0
\end{equation}
and consequently
\begin{equation*}
e^{\lambda t} S_A(t)x_0 = x_0 \ffor t\ge 0.
\end{equation*}
Writing $u_n(t) := u(t; s_n,\varepsilon_n,  x_n)$ for $n\ge 1$, by Theorem \ref{tw-con-21}, we assert that
\begin{equation}\label{c17}
u_n(t) \to u(t;0,0, x_0) \equiv x_0 \ \text{ uniformly for }t\in [0,T].
\end{equation}
On the other hand $$x_0\in\partial (U \oplus V) = \partial_{N_0} U \oplus V \cup U \oplus \partial_{X^\alpha_+\oplus X^\alpha_-} V,$$ and hence $x_0\in\partial_{N_0} U$ because \eqref{c11} implies that $x_0\in N_0$. Acting the operator $P$ on the equation
\begin{equation*}
x_n = e^{\lambda T}S_A(T)x_n + \varepsilon_n \int_0^T e^{\lambda (T - \tau)}S_A(T - \tau) G(s_n, \tau, u_n(\tau)) \d \tau,
\end{equation*}
and taking into account \eqref{ds1} and the inclusion $\ker(\lambda I - A)\subset \ker(I - e^{\lambda t}S_A(t))$
\begin{equation*}
\begin{aligned}
Px_n & = e^{\lambda T}S_A(T)Px_n + \varepsilon_n \int_0^T e^{\lambda (T - \tau)}S_A(T - \tau) P G(s_n, \tau, u_n(\tau)) \d \tau \\
& = Px_n + \varepsilon_n \int_0^T  P G(\tau, u_n(\tau)) \d \tau \ffor n\ge 1,
\end{aligned}
\end{equation*}
which implies that
\begin{equation*}
\int_0^T  P G(s,\tau, u_n(\tau)) \d \tau = 0 \ffor n\ge 1.
\end{equation*}
Letting $n\to+\infty$, by \eqref{c17}, we obtain
\begin{equation*}
h(x_0) = \int_0^T  P G(s,\tau, x_0) \d \tau = 0 \quad\text{ where } \ x_0\in \partial_{N_0} U,
\end{equation*}
which contradicts the assumption. Therefore there is $\ve_0 > 0$ such that for any $\varepsilon\in(0,\varepsilon_0]$, the map $\Psi_T(\varepsilon,\,\cdot \, , \, \cdot\,):[0,1]\times \o{U\oplus V} \to X^\alpha$ is an admissible homotopy and therefore
\begin{equation}\label{d3}
\begin{aligned}
\mathrm{deg_{LS}}(I - \Phi_T(\varepsilon,\,\cdot\,), U \oplus V) & = \degl(I - \Psi_T(\varepsilon,1, \, \cdot\,), U\oplus V) \\
& = \degl(I - \Psi_T(\varepsilon,0, \, \cdot\,), U\oplus V)
\end{aligned}
\end{equation}
for $\varepsilon\in(0,\varepsilon_0]$. Let $\phi^2_T:X^\alpha_+\oplus X^\alpha_-\to X^\alpha_+\oplus X^\alpha_-$ be an operator given by
\begin{equation*}
\phi^2_T(x) = e^{\lambda T}S_A(T)x \ffor x\in X^\alpha_+\oplus X^\alpha_-
\end{equation*}
and let $\phi^1_T(\ve, \cdot): N_0\to N_0$ the translation operator associated with
$$\dot u(t) = \ve h(u(t)), \qquad t > 0.$$
Then, it is not difficult to see that
$$\Psi_T(\varepsilon,0, x) = \phi^1_T(\ve, Px) + \phi^2_T(Q x) \ffor x\in X^\alpha,$$
and hence, for any $\ve\in(0,1]$, the map $\Psi_T(\varepsilon,0, \, \cdot\,):X^\alpha\to X^\alpha$ is topologically equivalent with $\t\Psi_T: [0,1]\times N_0\times (X^\alpha_+\oplus X^\alpha_-) \to N_0\times (X^\alpha_+\oplus X^\alpha_-)$ given by $$\t\Psi_T(\ve,u,v) = (\varphi^1_T(\ve, u), \varphi^2_T(v)) \ffor \ve\in [0,1], \ (u,v)\in N_0\times (X^\alpha_+\oplus X^\alpha_-)$$ and therefore
\begin{equation}\label{row111}
    \degl(I - \Psi_T(\varepsilon,0, \, \cdot\,), U\oplus V) = \degl(I - \t\Psi_T(\ve, \,\cdot\,), U\times V)
\end{equation}
for $\ve\in (0,\ve_0]$. Observe that Lemma \ref{lem-deg1} asserts that $\varphi^2_T(x) \neq x$ for $x\neq 0$ and
\begin{equation}\label{equ-deg}
\degl(I - \varphi^2_T, V) = \degl(I - e^{\lambda T}S_A(T)|_{X^\alpha_-\oplus X^\alpha_+}, V) = (-1)^{d_{k-1}}.
\end{equation}
Furthermore $h(x)\neq 0$ for $x\in\partial U$, and hence Theorem \ref{th-kras} says that there is $\ve_1 \in (0,\ve_0]$ such that $\varphi^1_T(\ve,x) \neq x$ for $\ve\in(0,\ve_1]$, $x\in \partial U$ and
\begin{equation}\label{row11}
    \degb(I - \varphi^1_T(\ve,\,\cdot\,), U) = \degb(-g, U).
\end{equation}
By \eqref{equ-deg}, \eqref{row11} and the multiplicative property of topological degree, for any $\ve \in (0,\ve_1]$, we have
\begin{equation*}
\begin{aligned}
\degl(I - \t\Psi_T(\ve, \,\cdot\,), U\times V) & = \degl (I - e^{\lambda T}S_A(T), V) \cdot \degb(I - \varphi^1_T(\ve,\,\cdot\,), U) \\
& \hspace{-12mm} = (-1)^{d_{k-1}} \degb(-g, U) = (-1)^{d_{k-1}} \cdot (-1)^{\dim N_0} \degb(g, U) \\
& \hspace{-12mm} = (-1)^{d_k} \degb(g, U).
\end{aligned}
\end{equation*}
Combining this with \eqref{d3}, \eqref{row111} we infer that
\begin{equation*}
\mathrm{deg_{LS}}(I - \Phi_T(\varepsilon,\,\cdot\,), U \oplus V) = (-1)^{d_k} \degb(g, U)
\end{equation*}
for $\ve \in (0,\ve_1]$, which completes the proof. \hfill $\square$ \\

An immediate consequence of Theorem \ref{th-aver-ker-sec} is the following corollary.

\begin{corollary}\label{cor-aver}
Let $U \subset N_0$ and $V \subset X^\alpha_+\oplus X^\alpha_-$ where $0\in V$, be open bounded sets such that $h(x)\neq 0$ for $x \in \partial_{N_0} U$. If $\mathrm{deg_B}(g, U) \neq 0$ then there is $\varepsilon_0 \in (0,1)$ such that, for any $\varepsilon\in (0,\varepsilon_0]$, the equation \eqref{równanie-gl-eps} admits a $T$-periodic mild solution.
\end{corollary}

\section{Index formula for periodic solutions}

We will study the problem of existence of $T$-periodic solutions for the equation
\begin{equation}\label{row-tt}
\dot u(t) = - A u(t) + \lambda u(t) + F (t,u(t)), \qquad  t > 0,
\end{equation}
where $\lambda$ is an eigenvalue of the operator $A:X\supset D(A) \to X$ and $F\colon X^\alpha \to X$ is a continuous map. Assume that assumptions $(A1)$, $(A2)$, $(A3)$, $(F1)$, $(F3)$, $(F4)$ hold and furthermore \\[5pt] \noindent\makebox[9mm][l]{$(F5)$}\parbox[t][][t]{118mm}{$F(t + T,x) = F(t,x)$ for $t\in [0,+\infty)$ and $x\in X^\alpha$.} \\[5pt]
Theorem \ref{th-exist1} implies that, for any $x\in X^\alpha$, there is mild solution $u(\,\cdot\,;x):\R\to X^\alpha$ of \eqref{row-tt} starting at $x$. Let $\Phi_T: X^\alpha \to X^\alpha$ be an associated translation along trajectories operator give by $$\Phi_T(x) := u(T;x) \ffor x\in X^\alpha.$$ Then Theorems \ref{tw-con-21} and \ref{tw-exi-con-comp} say that $\Phi_T$ is a completely continuous map.

We say that solution $u:[0,+\infty)\to X^\alpha$ of the equation \eqref{row-tt} is $T$-periodic provided $u(t) = u(t + T)$ for $t\ge 0$. It is not difficult to check that every fixed point of the translation operator $\Phi_T$ can be identified with a starting point of $T$-periodic solution of \eqref{row-tt}.
\begin{remark}\label{rem-non-ex2}
If the equation \eqref{row-tt} is at resonance at infinity, the problem of existence of $T$-periodic solution may not have solutions for general nonlinearity $F$. \\[5pt]
To see this take $F(t,x) := y_0$ for $t\in[0,+\infty)$, $x\in X^\alpha$, where $y_0\in\ker(\lambda I - A)\setminus\{0\}$. If $u:[0,+\infty)\to X^\alpha$ is a $T$-periodic solution of \eqref{row-tt}, then we have the integral formula $$u(t) = e^{\lambda t}S_A(t)u(0) + \int_0^t e^{\lambda (t - \tau)}S_A(t - \tau)y_0\d\tau\ffor t\ge 0.$$ Since $\ker(\lambda I - A) \subset \ker(I - e^{\lambda t}S_A(t))$ for $t\ge 0$ we infer that $u(T) = e^{\lambda T}S_A(T)u(0) + T y_0$. Acting on this equation by the operator $P$ and using \eqref{ds1} we deduce that $$Pu(T) = e^{\lambda T} S_A(T)Pu(0) + T y_0  = Pu(0) + T y_0 = Pu(T) + T y_0,$$ which is impossible, because $y_0\neq 0$. \hfill $\square$
\end{remark}
To overcome difficulties described in the above remark we formulate the following \emph{geometrical conditions}, which will be used to prove that existence of $T$-periodic solutions for \eqref{row-tt}: \label{g1g2} \\[1pt]
\begin{equation*}\leqno{(G1)}
\quad\left\{\begin{aligned}
& \text{given set } B\subset X^\alpha_+\oplus X^\alpha_- \text{ there exists } R > 0 \text{ such that} \\
& \<F(t, x + y), x\> > 0 \text{ for } (t,y, x)\in [0,T]\times B \times X_0 \ \text{ where } \|x\|_H\ge R,
\end{aligned}\right.
\end{equation*}
\begin{equation*}\leqno{(G2)}
\quad\left\{\begin{aligned}
& \text{given set $B\subset X^\alpha_+\oplus X^\alpha_-$ there exists $R > 0$ such that} \\
& \<F(t, x + y), x\> < 0 \text{ for } (t,y, x)\in [0,T]\times B \times X_0 \text{ where } \|x\|_H\ge R.
\end{aligned}\right. 
\end{equation*}
Now we are ready to prove the following \emph{index formula for periodic solutions}, which determines the Leray-Schauder topological degree of the vector field $I - \Phi_T$ with resect to the ball with sufficiently large radius. This theorem is a tool to searching the fixed points of $\Phi_T$ and hence the $T$-periodic solutions for \eqref{row-tt}.
\begin{theorem}\label{th-reso-m}
Assume that $\lambda = \lambda_k$ for some $k\ge 1$.
\begin{enumerate}
\item[(i)] If condition $(G1)$ holds, then there is $R > 0$ such that $\Phi_T(x)\neq x$ for $x\in X^\alpha$ with $\|x\|_\alpha \ge R$ and
\begin{equation*}
    \mathrm{deg_{LS}}(I - \Phi_T, B(0,R)) = (-1)^{d_k}.
\end{equation*}
\item[(ii)] If condition $(G2)$ holds, then there is $R > 0$ such that $\Phi_T(x)\neq x$ for $x\in X^\alpha$ with $\|x\|_\alpha \ge R$ and
\begin{equation*}
    \mathrm{deg_{LS}}(I - \Phi_T, B(0,R)) = (-1)^{d_{k-1}}.
\end{equation*}
\end{enumerate}
Here $d_0 := 0$ and $d_l := \sum_{i=1}^l \dim\ker(\lambda_i I - A)$ if $l\ge 1$.
\end{theorem}
In the proof of this theorem we will consider the family of differential equations
\begin{equation}\label{diff:2}
\dot u(t) = - A u(t) + \lambda u(t) + \varepsilon F(t, u(t)), \qquad  t > 0
\end{equation}
where $\ve\in[0,1]$ is a parameter. Let $\Psi_T:[0,1] \times X^\alpha\to X^\alpha$ given by $$\Psi_T(\ve, x) := u(T;\ve, x) \ffor \ve\in[0,1], \ x\in X^\alpha,$$ be the translation along trajectories operator associated with this equation. \\[5pt]
Before we proceed to prove of the above theorem we prove the following lemmata which provides {\em a priori} estimates on $T$-periodic mild solutions.
\begin{lemma}\label{lem-boun}
There is a constant $R > 0$ such that if $u:=u_\varepsilon:[0,+\infty)\to X^\alpha$, where $\varepsilon\in(0,1]$, is a $T$-periodic mild solution for \eqref{diff:2}, then
\begin{equation}\label{nier11}
    \|Qu(t)\|_\alpha \le R \ffor t\in [0,T].
\end{equation}
\end{lemma}
\noindent\textbf{Proof.} Since $u$ is $T$-periodic, for any integer $k > 0$, we have equality
\begin{equation*}
    u(t) = u(t + kT) \ffor t\in [0,T],
\end{equation*}
which implies that
\begin{equation}\label{d-row}
u(t) = e^{\lambda k T}S_A(kT)u(t) + \varepsilon\int_t^{t + kT} e^{\lambda (t + kT - \tau)}S_A(t + kT - \tau)F(\tau, w(\tau)) \d \tau
\end{equation}
for $t\ge 0$ and $k\ge 1$. Acting on \eqref{d-row} by the operator $Q_+$ and using \eqref{ds1}, we obtain
\begin{equation*}
Q_+u(t) = e^{\lambda t}S_A(t)Q_+ u(t) + \varepsilon\int_t^{t+T} e^{\lambda (t + T - \tau)}S_A(t + T - \tau) Q_+ F(\tau, u(\tau)) \d \tau
\end{equation*}
for $t\ge 0$ and $n\ge 1$. If $m$ is the constant from assumption $(F3)$, then, by \eqref{ine11}
\begin{equation*}
\begin{aligned}
\|Q_+u(t)\|_\alpha & \le \|e^{\lambda k T}S_A(kT)Q_+u(t)\|_\alpha \\
& \hspace{5mm} + \int_t^{t + k T} \|A_\delta^\alpha e^{\lambda(t + k T - \tau)}S_A(t + k T - \tau)Q_+F(\mu,u(\tau))\| \d \tau \\
& \le \|e^{\lambda k T}S_A(k T)Q_+u(t)\|_\alpha + M \int_t^{t + k T} \frac{e^{- c (t + k T - \tau)}}{(t + k T - \tau)^\alpha} \, \|Q_+F(\mu,u(\tau))\| \d \tau \\
& \le \|e^{\lambda k T }S_A(k T)Q_+u(t)\|_\alpha + m M \|Q_+\|_{L(X)} \int_t^{t + k T} \frac{e^{- c (t + k T - \tau)}}{(t + k T - \tau)^\alpha} \d \tau \\
& \le M \frac{e^{- c kT }}{(kT)^\alpha} \, \|Q_+u(t)\| + m M \|Q_+\|_{L(X)} \int_t^{t + k T} \frac{e^{- c (t + k T - \tau)}}{(t + k T - \tau)^\alpha} \d \tau.
\end{aligned}
\end{equation*}
where $c,M > 0$. On the other hand, for $k \ge 2$, we have
\begin{equation*}
\begin{aligned}
\int_t^{t + kT} \frac{e^{- c (t + kT - \tau)} }{(t + kT - \tau)^\alpha} \d \tau
& = \int_t^{t + (k-1)T} \frac{e^{- c (t + kT - \tau)}}{(t + kT - \tau)^\alpha} \d \tau
  + \int_{t + (k-1)T}^{t + kT} \frac{e^{- c (t + kT - \tau)}}{(t + kT - \tau)^\alpha} \d \tau \\
& \hspace{-5mm} \le \int_t^{t + (k-1)T} T^{-\alpha} e^{- c (t + kT - \tau)} \d \tau
 + \int_{t + (k-1)T}^{t + kT} \frac{1}{(t + kT - \tau)^\alpha}  \d \tau \\
& \hspace{-5mm} = T^{-\alpha}(e^{-c T} - e^{-c kT})/c + T^{1-\alpha}/(1 - \alpha).
\end{aligned}
\end{equation*}
In a consequence, for any $t\in [0,T]$ and integer $k > 0$, we find that
\begin{equation*}
\|Q_+u(t)\|_\alpha \le M \frac{e^{- c kT }}{(kT)^\alpha} \, \|Q_+u(t)\| + mM \|Q_+\|_{L(X)}T^{-\alpha}\left((e^{-c T} - e^{-c k T})/c + \frac{T}{1 - \alpha}\right).
\end{equation*}
Hence, letting with $k\to +\infty$, we assert that
\begin{equation}\label{eq:31}
\|Q_+u(t)\|_\alpha \le m M \|Q_+\|_{L(X)} T^{-\alpha} \left(  e^{-c T}/c + \frac{T}{1 - \alpha}\right)=:R_1
\end{equation}
for $t\in [0,T]$. Now we act on equation \eqref{d-row} by operator $Q_-$. Then, in view of \eqref{ds1}, we have
\begin{equation}\label{eq:11}
\hspace{-1mm} e^{-\lambda k T}S_A(-k T) Q_-u(t) = Q_-u(t) + \int_t^{t + kT} e^{\lambda (t - \tau)},
S_A(t - \tau)Q_-F(\mu,u(\tau)) \d \tau
\end{equation}
for any $t\in [0,1]$ and integer $k\ge 1$, because the semigroup $\{S_A(t)\}_{t\ge 0}$ extend on $X_-$ to a $C_0$ group of bounded operators. Therefore, by \eqref{ine33}
\begin{align*}
\|e^{-\lambda k T}S_A(-k T) Q_-u(t)\| \le M \, e^{-c k T} \|Q_-u(t)\|,
\end{align*}
where $c,M > 0$, which together with \eqref{eq:11} gives
\begin{align*}
\|Q_-u(t)\| & \le \|e^{-\lambda k T}S_A(- k T) Q_-u(t)\| \\
& \qquad + \int_t^{t + k T} \|e^{\lambda (t - \tau)} S_A(t - \tau)Q_-F(s,u(\tau))\| \d \tau \\
& \le  M \, e^{-c k T} \|Q_-u(t)\| + M \int_t^{t + kT} e^{c (t - \tau)}\|Q_-F(s,u(\tau))\| \d \tau \\
& \le  M \, e^{-c k T} \|Q_-u(t)\| + m M \|Q_-\|_{L(X)} \int_t^{t + kT}  e^{c (t - \tau)} \d \tau \\
& =  M \, e^{-c k T} \|Q_-u(t)\| + m M \|Q_-\|_{L(X)}\left(1 - e^{ - c k T}\right)/c.
\end{align*}
Letting $k\to +\infty$ we obtain
\begin{equation}\label{eq:32}
\|Q_-u(t)\| \le m M \|Q_-\|_{L(X)}/c \ffor t\in [0,T].
\end{equation}
Since $X_-$ is finite dimensional, there is a constant $C' > 0$ such that
\begin{equation}\label{norm-in22}
    \|x\|_\alpha \le C'\|x\| \ffor x\in X_-.
\end{equation}
Hence, by \eqref{eq:32}, we have
\begin{equation}\label{eq:32bbb}
\|Q_-u(t)\|_\alpha \le m C' M \|Q_-\|_{L(X)}/c =:R_2 \ffor t\in [0,T].
\end{equation}
Taking into account the inequalities \eqref{eq:31} and \eqref{eq:32bbb}, for any $t\in [0,T]$, we have
\begin{align*}
\|Qu(t)\|_\alpha & \le \|Q_-u(t)\|_\alpha + \|Q_+u(t)\|_\alpha \le \|Q_-\|_{L(X^\alpha)}\|u(t)\|_\alpha + \|Q_+\|_{L(X^\alpha)}\|u(t)\|_\alpha  \\ & \le R_1 \|Q_+\|_{L(X^\alpha)} + R_2 \|Q_-\|_{L(X^\alpha)}:= R,
\end{align*}
which completes the proof. \hfill $\square$ \\

\begin{lemma}\label{th-guid-fun}
Assume that $N_0:=\ker(A - \lambda I)$ and let $g\colon N_0 \to N_0$ be given by
\begin{equation*}
g(x) := \int_0^T PF(s,x)\,d s \qquad\mathrm{for}\quad x\in N_0.
\end{equation*}
\begin{enumerate}
\item[(i)] If condition $(G1)$ holds, then there is $R_0 > 0$ such that $g(x)\neq 0$ for $x\in N_0$ with $\|x\|_H \ge R_0$ and
\begin{equation*}
\mathrm{deg_B}(g, B(0,R)) = 1 \qquad\mathrm{for}\quad R\ge R_0.
\end{equation*}
\item[(ii)] If condition $(G2)$ holds, then there is $R_0 > 0$ such that $g(x)\neq 0$ for $x\in N_0 $ with $\|x\|_H \ge R_0$ and
\begin{equation*}
\mathrm{deg_B}(g, B(0,R)) = (-1)^{\dim N_0 } \qquad\mathrm{for}\quad R\ge R_0.
\end{equation*}
\end{enumerate}
\end{lemma}
\noindent\textbf{Proof.} For the proof of $(i)$, define the map $H\colon [0,1]\times N_0 \to N_0 $ by \[H(s,x) := s g(x) + (1 - s)x \qquad\text{for}\quad x\in N_0.\] By condition $(G1)$ there is a constant $R_0 > 0$ such that
\begin{equation*}
\<F(\tau, x), x\>_H  > 0 \quad\text{ for } \ \ \tau\in [0,T], \  x\in N_0 \ \text{ such that } \|x\|_H\ge R_0,
\end{equation*}
which, after integration, implies that
\begin{equation}\label{wa1a}
    \<g(x),x\>_H  = \int_0^T \<F(\tau, x), x\>_H \d\tau > 0 \ \ \text{ for } \ x\in N_0 \ \text{ such that } \|x\|_H\ge R_0.
\end{equation}
Let $R \ge R_0$. We show that $H(s,x)\neq 0$ for $s\in [0,1]$ and $x\in N_0 $ where $\|x\|_H = R$. Otherwise there is $s\in [0,1]$ and $x\in N_0$ with $\|x\|_H = R$ such that $H(s,x) = 0$. Consequently
\begin{align*}
0 = \langle H(s,x),x \rangle_H = s \langle g(x),x \rangle_H + (1 - s) \langle x,x \rangle_H.
\end{align*}
If $s = 0$ then $0=\|x\|_H^2 = R^2$, which is impossible. If $s\in (0,1]$ then
$0 \ge \langle g(x),x \rangle$, which contradicts \eqref{wa1a}. Hence, by the homotopy invariance,
\begin{align*}
\mathrm{deg_B}(g, B(0,R)) & = \mathrm{deg_B}(H(1,\,\cdot\,), B(0,R)) = \mathrm{deg_B}(H(0,\,\cdot\,), B(0,R)) \\
& = \mathrm{deg_B}(I, B(0,R)) = 1,
\end{align*}
and the proof of $(i)$ is completed. To verify $(ii)$ observe that condition $(G2)$ implies the existence of $R_0 > 0$ such that
\begin{equation*}
\<F(\tau, x), x\>_H  < 0 \quad\text{ for } \ \  \tau\in [0,T], \ x\in N_0 \ \text{ with } \|x\|_H\ge R_0,
\end{equation*}
which, after integration, gives
\begin{equation}\label{wa1b}
    \<g(x),x\>_H  = \int_0^T\<F(\tau, x), x\>_H \d\tau< 0 \ffor x\in N_0 \ \text{ with } \|x\|_H\ge R_0.
\end{equation}
Therefore, for any $R > R_0$, the homotopy $H\colon [0,1]\times N_0  \to N_0 $ given by
\[H(s,x) := s g(x) - (1 - s)x \qquad\text{for}\quad x\in N_0 \] is such that $H(s,x) \neq 0$ for $s\in [0,1]$ and $x\in N_0$ with $\|x\|_H = R$. Indeed, if $H(s,x) = 0$ for some $s\in [0,1]$ and $x\in N_0 $ with $\|x\|_H = R$, then
\begin{align*}
0 = \langle H(s,x),x \rangle_H = s \langle g(x),x \rangle_H - (1 - s)\langle x,x \rangle_H.
\end{align*}
If $s\in (0,1]$ then $\langle g(x),x \rangle_H \ge 0$, contrary to \eqref{wa1b}. If $s=0$ then
$R^2 = \|x\|_H^2 = 0$, which is again impossible. Hence, by the homotopy invariance,
\begin{equation*}
\mathrm{deg_B}(g, B(0,R)) = \mathrm{deg_B}(-I, B(0,R)) = (-1)^{\dim N_0 },
\end{equation*}
which completes the proof. \hfill $\square$\\

\noindent\textbf{Proof of Theorem \ref{th-reso-m}.} \textbf{Step 1.} We show that there is $R_0 > 0$ such that
\begin{equation}\label{co:7}
    \Psi_T(\ve, x) \neq x \ \ \text{for} \ \ \ve\in (0,1] \ \text{ and } \ x \in X^\alpha \ \text{ with } \
    \|x\|_\alpha \ge R_0.
\end{equation}
Suppose contrary that there are sequences $(x_n)$ in $X^\alpha$ and $(\ve_n)$ in $(0,1]$ such that $\|x_n\|_\alpha \to +\infty$ as $n\to +\infty$ and
\begin{equation*}
\Psi_T(\ve_n, x_n) = x_n \ffor n\ge 1.
\end{equation*}
Writing $z_n:= x_n/\|x_n\|_\alpha$, $u_n := u(\,\cdot\,;\ve_n, x_n)$ and $v_n := u(\,\cdot\,;\ve_n, x_n)/\|x_n\|_\alpha$ we see that
\begin{equation}\label{eq:1}
v_n(t) = e^{\lambda t}S_A(t)z_n + \varepsilon_n\int_0^t e^{\lambda (t - \tau)}S_A(t - \tau)F(\tau, u_n(\tau))/\|u_n\|_\alpha \d \tau
\end{equation}
for $t\in[0,T]$ and $n\ge 1$. Define
\begin{equation}\label{eq:22}
y_n(t):=\varepsilon_n\int_0^t e^{\lambda (t - \tau)}S_A(t - \tau)F(\tau, u_n(\tau))/\|x_n\|_\alpha \d \tau.
\end{equation}
Since $A$ is sectorial, there are constants $M > 0$ and $c_0\in \R$ such that $$\|A_\delta^\alpha S_A(t)x\| \le Mt^{-\alpha} e^{c_0 t} \ffor t > 0.$$ Then, for any $t\in[0,T]$ and $n\ge 1$, we have
\begin{align*}
\|y_n(t)\|_\alpha & \le \varepsilon_n\int_0^t \|e^{\lambda (t - \tau)} A_\delta^\alpha S_A(t - \tau)F(\tau, u_n(\tau))\|/\|x_n\|_\alpha \d \tau \\
& \le \int_0^t M e^{(|\lambda| + |c|) T}(t - \tau)^{-\alpha}\|F(\tau, u_n(\tau))\|/\|x_n\|_\alpha \d \tau \\
& \le \int_0^t m M e^{(|\lambda| + |c|) T}(t - \tau)^{-\alpha}/\|x_n\|_\alpha \d \tau
\le  \frac{m Me^{(|\lambda| + |c|) T}}{(1 - \alpha)\|x_n\|_\alpha} \, T^{1 - \alpha},
\end{align*}
where $m$ is the constant from $(F3)$. Therefore, letting $n\to +\infty$, we find that
\begin{equation}\label{eq:23}
  \|y_n(t)\|_\alpha \to 0 \qquad\text{ as }\quad n\to +\infty,
\end{equation}
and the convergence in uniform for $t\in [0,T]$. With the compactness of the semigroup $\{S_A(t)\}_{t\ge 0}$ the set $\{e^{\lambda t}S_A(t)z_n \ | \ n\ge 1\}$ is relatively compact in $X^\alpha$. Combining this fact with \eqref{eq:1} and \eqref{eq:23}, we infer that the set $$\{z_n \ | \ n\ge 1\} = \{v_n(T) \ | \ n\ge 1\}$$ is also relatively compact in $X^\alpha$. Therefore, passing if necessary to a subsequence, we can assume that there is $z_0 \in X^\alpha$ such that $z_n \to z_0$ in $X^\alpha$ as $n\to +\infty$. Since $\|z_n\|_\alpha = 1$, we obtain $\|z_0\|_\alpha = 1$. Writing the equality \eqref{eq:1} with $t := T$ we have
\begin{equation*}
z_0 = e^{\lambda T}S_A(T)z_0
\end{equation*}
which, by \eqref{row-sp}, implies that $z_0 \in\ker(\lambda I - A)$ and consequently
\begin{equation*}
e^{\lambda t}S_A(t)z_0 = z_0 \ffor t\ge 0.
\end{equation*}
From \eqref{eq:23}, we conclude that, for $t\in [0,T]$,
\begin{equation}\label{eq:25}
v_n(t) \to z_0 \ \text{ w }X^\alpha, \text{ as } n \to +\infty \ \text{ uniformly for } \ t\in[0,T].
\end{equation}
From Lemma \ref{lem-boun} we deduce that there is a constant $C > 0$ such that
\begin{equation}\label{imp-ine}
    \|Qu_n(t)\| \le C \ffor t\in [0,T], \ n\ge 1.
\end{equation}
In conditions $(G1)$ and $(G2)$, let the set $B$ be a ball in $X^\alpha$ with the radius $C$. Using \eqref{ort}, we infer that there is $R_0 > 0$ such that
\begin{equation}\label{ro2}
\<PF(t, x + y), x\> > 0 \quad\text{ for } \ \  (t,y, x)\in [0,T]\times B \times X_0 \ \text{ with } \|x\|_H\ge R_0,
\end{equation}
if condition $(G1)$ is satisfied and
\begin{equation}\label{ro2bb}
\<PF(t, x+y), x\> < 0 \quad\text{ for } \ \  (t,y, x)\in [0,T]\times B \times X_0 \ \text{ with } \|x\|_H\ge R_0,
\end{equation}
if the condition $(G2)$ holds. Acting by the operator $P$ on the equation
\begin{equation*}
u_n(t) = e^{\lambda t}S_A(t) x_n + \varepsilon_n\int_0^t e^{\lambda (t - \tau)}S_A(t - \tau)F(\tau, u_n(\tau))\d \tau \ffor t\ge 0,
\end{equation*}
from \eqref{ds1} and the inclusion $\ker(\lambda I - A)\subset \ker(I - e^{\lambda t}S_A(t))$, we have
\begin{equation*}
Pu_n(t) =  P x_n + \varepsilon_n\int_0^t  P F(\tau, u_n(\tau)) \d \tau
\end{equation*}
for $t\in[0,T]$ and $n\ge 1$. Hence the map $Pu_n$ is continuously differentiable on $[0,T]$ and
\begin{equation*}
\frac{d u_n(t)}{dt} = \varepsilon_n P F(t, u_n(t)) \ffor n\ge 1.
\end{equation*}
Therefore, for $t\in[0,T]$ and $n\ge 1$, we have
\begin{equation*}
\frac{d}{dt}\frac{1}{2} \|P u_n(t)\|^2_H = \left\<\frac{d u_n(t)}{dt}, u_n(t)\right\>_H = \varepsilon_n \< P F(t, u_n(t)), P u_n(t)\>_H
\end{equation*}
which, after integration, gives
\begin{equation}\label{eq:18a}
\begin{aligned}
0 & = \frac{1}{2}(\|Pu_n(T)\|_H - \|Pu_n(0)\|_H) = \varepsilon_n \int_0^T \< P F(\tau, u_n(\tau)), Pu_n(\tau)\>_H \d \tau \\
& = \varepsilon_n \int_0^T \< F(\tau, Qu_n(\tau) + \|u_n\|_\alpha Pv_n(\tau)), \|u_n\|_\alpha Pv_n(\tau)\>_H \d \tau
\end{aligned}
\end{equation}
for $n\ge 1$. In view of \eqref{eq:25} we find that $Pv_n(t) \to Pz_0 = z_0$ in $X^\alpha$, uniformly for $t\in[0,T]$. Since  $z_0\neq 0$, there is $n_0\ge 1$ such that $\|Pv_n(t) - z_0\|_H\le \|z_0\|_H/2$ for $n\ge n_0$ and $t\in[0,T]$. Then $$\|Pv_n(t)\|_H \ge \|z_0\|_H - \|z_0\|_H/2 = \|z_0\|_H/2 \ffor n\ge n_0, \ t\in[0,T],$$ and hence, increasing $n_0 \ge 1$ if necessary, we deduce that
\begin{equation}\label{rowna11}
\|\|u_n\|_\alpha Pv_n(t)\| \ge R_0 \ffor n\ge n_0, \ t\in[0,T].
\end{equation}
In the case of point $(i)$, the inequality \eqref{rowna11} together with \eqref{ro2} and \eqref{imp-ine}, imply that
\begin{align*}
\int_0^T \< PF(\tau, Qu_n(\tau) + \|u_n\|_\alpha Pv_n(\tau)), \|u_n\|_\alpha Pv_n(\tau)\> \d \tau  > 0 \ffor n\ge n_0.
\end{align*}
On the other hand, in the case of point $(ii)$, the inequality \eqref{rowna11} along with \eqref{ro2bb} and \eqref{imp-ine}, give
\begin{align*}
\int_0^T \< PF(\tau, Qu_n(\tau) + \|u_n\|_\alpha Pv_n(\tau)), \|u_n\|_\alpha Pv_n(\tau)\> \d \tau  < 0 \ffor n\ge n_0.
\end{align*}
In the both cases we obtain a contradiction with \eqref{eq:18a}, because $\ve_n\in(0,1]$ for $n\ge 1$. Thus we proved  \eqref{co:7} and the proof of Step 1 is completed. \\[5pt]
\textbf{Step 2.} To complete the proof of theorem, we show that there is $\ve_0 > 0$ such that for any $\ve\in (0,\ve_0]$
\begin{equation}\label{ro-de1}
\mathrm{deg_{LS}}(I - \Psi_T(\ve,\,\cdot\,), B(0,R_0)) = (-1)^{d_k}
\end{equation}
if condition $(G1)$ is satisfied and
\begin{equation}\label{ro-de2}
\mathrm{deg_{LS}}(I - \Psi_T(\ve,\,\cdot\,), B(0,R_0)) = (-1)^{d_{k-1}} ,
\end{equation}
if condition $(G2)$ holds. Lemma \ref{th-guid-fun} asserts the existence of $R_1 > R_0$ such that
\begin{equation}\label{e1}
g(x) \neq 0 \ \ \text{for} \ \ x\in N_0 \ \text{ such that } \ \|x\|_H \ge R_1
\end{equation}
and furthermore
\begin{equation}\label{w-G5} 
\mathrm{deg_B}(g, B(0,R)) = 1 \qquad\mathrm{for}\quad R\ge R_1,
\end{equation}
if condition $(G1)$ is satisfied and
\begin{equation}\label{w-co6}
\mathrm{deg_B}(g, B(0,R)) = (-1)^{\dim N} \qquad\mathrm{for}\quad R\ge R_1
\end{equation}
if condition $(G2)$ holds. Let $R_2:=\max(R_1/C_1, R_1)$, where $C_1>0$ is a constant such that
\begin{equation}\label{wzor1}
C_1(\|P x\|_H + \|Q x\|_\alpha) \le \|x\|_\alpha \ffor x\in X^\alpha.
\end{equation}
Define $$U:=\{x\in N_0 \ | \ \|x\|_H \le R_2\} \ \ \text{ and } \ \ V:=\{x\in X^\alpha_+\oplus X^\alpha_- \ | \ \|x\|_\alpha \le R_2\}.$$ In view of \eqref{wzor1} we deduce that $B(0,R_1)\subset U\oplus V$. By Step 1 and the fact that $R_1 > R_0$, the excision property of topological degree gives
\begin{equation}\label{de2}
\degl(I - \Psi_T(\ve,\,\cdot\,), B(0,R_0)) = \degl(I - \Psi_T(\ve,\,\cdot\,), U\oplus V)
\end{equation}
for $\ve\in(0,1]$. Further, from \eqref{e1} and the fact that $R_2 \ge R_1$ we find that $g(x)\neq 0$ for $x\in \partial_{N_0} U$. Hence, by Theorem \ref{th-aver-ker-sec} we have $\ve_0 \in (0,1)$ such that for any $\ve\in (0,\ve_0]$, $\Psi_T(\ve,x)\neq x$ for $x\in \partial(U \oplus V)$ and
\begin{equation*}
\mathrm{deg_{LS}}(I - \Psi_T(\ve,\,\cdot\,), U \oplus V) = (-1)^{d_k} \cdot \mathrm{deg_B}(g,U),
\end{equation*}
which, by \eqref{de2}, gives
\begin{equation*}
\degl(I - \Psi_T(\ve,\,\cdot\,), B(0,R_0)) = (-1)^{d_k} \cdot \mathrm{deg_B}(g,U) \ffor \ve\in(0,\ve_0].
\end{equation*}
Combining this with \eqref{w-G5} and \eqref{w-co6}, we prove \eqref{ro-de1} and \eqref{ro-de2}, which completes the proof. \hfill $\square$ \\

\section{Applications} 
 
Let us assume that $\Omega\subset\R^n$ is an open bounded set with $C^\infty$ boundary. Let $\mathcal{A}$ be a uniformly elliptic symmetric second order differential operator with a Dirichlet boundary conditions:
$$\mathcal{A} \y (x) = - \sum_{i,j=1}^n D_j(a_{ij}(x)D_i \y(x)) \ffor \y\in C^2(\o\Omega)$$ 
with $a_{ij} = a_{ji}\in C^2(\overline\Omega)$ for $1\le i,j\le n$.  
Furthermore let $g:[0,+\infty)\times\Omega\times\R\times\R^n\to\R$ be a continuous map satisfying the following assumptions: \\[5pt]
\noindent\makebox[22pt][l]{$(E1)$} \parbox[t][][t]{143mm}{there exists $C > 0$ such that 
    \begin{align*}
    |g(t,x,s_1,y_1) - g(t,x,s_2,y_2)| & \le C(|s_1 - s_2| + |y_1 - y_2|),
    \end{align*}
    for $t\in [0,+\infty)$, $x\in\Omega$, $s_1,s_2\in\R$ and $y_1,y_2\in\mathbb{R}^n$,}\\[3pt] 
\noindent\makebox[22pt][l]{$(E2)$} \parbox[t][][t]{143mm}{there is a constant $m > 0$ such that 
$$|g(t,x,s,y)| \le m\ffor x\in\Omega, \ \in\R, \ y\in\mathbb{R}^n, \ t\in [0,+\infty).$$}\\[5pt]
Write $X:=L^p(\Omega)$, where $p\ge 1$, and define the operator $A_p: X\supset D(A_p)\to X$ by
\begin{equation}
\begin{aligned}\label{op-a}
D(A_p) := W^{2,p}_0(\Omega) \quad\text{and}\quad A_p \y  := \mathcal{A} \y \ffor \y\in D(A_p).
\end{aligned}
\end{equation}

\begin{proposition}{\em (see \cite{MR1778284, MR500580})}\label{th-sect} The following assertions hold. \\[5pt]
\noindent\makebox[5mm][l]{$(a)$} \parbox[t][][t]{143mm}{The operator $A_p$ is positively defined, sectorial and has compact resolvent.}\\[5pt]
\noindent\makebox[5mm][l]{$(b)$} \parbox[t][][t]{143mm}{If the domain $D(A_p)$ is equipped with the graph norm $$\|\y\|_{D(A_p)}:= \|A_p\y\|_{L^p(\Omega)} + \|\y\|_{L^p(\Omega)} \ffor \y\in D(A_p),$$ then inclusion $D(A_p) \subset W^{2,p}(\Omega)$ is compact.}\\[5pt]
\noindent\makebox[5mm][l]{$(c)$} \parbox[t][][t]{143mm}{The operator $A_2 :L^2(\Omega) \supset D(A_2) \to L^2(\Omega)$ is self-adjoint.}
\end{proposition}
From Lemma \ref{th-sect} $(a)$ it follows that the operator $A_p:X\supset D(A_p)\to X$ is positively defined sectorial operator on $X$, and hence it define a fractional space $X^\alpha := D(A_p^\alpha)$, ($\alpha\in(0,1)$) with the norm
\begin{align*}
\|\y\|_\alpha := \|A_p^\alpha \y\| \ffor \y\in X^\alpha.
\end{align*}
From now on we assume that \\[5pt]
\noindent\makebox[22pt][l]{$(E4)$} \parbox[t][][t]{143mm}{$p\ge 2n$ and $\alpha\in(3/4,1)$.}
\begin{remark}\label{rem-pom2} 
$(a)$ Observe that $A_p$ satisfies assumptions $(A1)$, $(A2)$ and $(A3)$. Indeed, by Theorem \ref{th-sect} $(a)$ we infer that $A_p$ has compact resolvent, that is, $(A1)$ holds. Let us take $H:=L^2(\Omega)$ equipped with the standard inner product and norm. Since $\Omega$ is bounded and $p\ge 2$ we have the embedding $i:L^p(\Omega) \hookrightarrow L^2(\Omega)$ and hence assumption $(A2)$ is satisfied. Using the boundedness of $\Omega$ again, we see that for $\h A:= A_2$ we have $D(A_p)\subset D(\h A)$ and furthermore $\h A \y = A_p \y$ for $\y \in D(A_p)$. This proves that $i\times i\left[\mathrm{Gr}\,A_p \right]\subset \h A$. By Theorem \ref{th-sect} $(c)$ the operator $\h A$ is self-adjoint and therefore the assumption $(A3)$ is satisfied. \\[5pt] 
$(b)$ By Remark \ref{rem-pom}, we see that $\sigma(A_p)= \{\lambda_i\}$ where $$0 < \lambda_1 < \lambda_2 < \ldots < \lambda_i < \lambda_{i+1} < \ldots$$ is a sequence of eigenvalues, which is finite or $\lambda_i \to +\infty$ as $i\to +\infty$. \\[5pt]
$(c)$ In view of $(E4)$ one has $\alpha\in (3/4, 1)$ and $p\ge 2n$. Therefore $2\alpha - \frac{n}{p} > 1$ and using \cite[Theorem 1.6.1]{MR610244} one has 
\begin{equation}\label{zan}
    X^\alpha\subset C^1(\o\Omega).
\end{equation}    
$(d)$ By \cite[Theorem 1.4.8]{MR610244}, the inclusion $X^\alpha \subset X^\beta$ is compact, if $\alpha > \beta > 0$.   \hfill $\square$
\end{remark}

By Remark \ref{rem-pom2} $(c)$, we can introduce a mapping $F\colon [0,+\infty)\times X^\alpha \to X$ given, for $\y\in X^\alpha$ by the following formula 
\begin{equation}\label{odwz-f}
    F(t,\y)(x) := g(t,x, \y(x), \nabla \y(x)) \ffor t\in [0,+\infty), \ x\in\Omega.
\end{equation}
We call $F$ \emph{the Niemytzki operator} associated with $f$.
\begin{lemma}\label{lem-nem-op}
The following assertions hold.
\begin{enumerate}
\item[(i)] The map $F$ is well defined, continuous and satisfies assumption $(F1)$. \\[-9pt]
\item[(ii)] There is $K > 0$ such that
\begin{equation}\label{sw}
\|F(t,\y)\| \le K \ffor t\in[0,+\infty), \ \y\in X^\alpha.
\end{equation}
\end{enumerate}
\end{lemma}
\noindent\textbf{Proof.} In view of \eqref{zan}, the inclusion $X^\alpha\subset W^{1,p}(\Omega)$ is continuous and hence there is $M > 0$ such that
\begin{equation*}
    \|\y\|_{W^{1,p}(\Omega)} \le M \|\y\|_\alpha \ffor \y\in X^\alpha.
\end{equation*}
By the assumption $(E2)$, for any $u\in X^\alpha$ we have
\begin{equation}\label{og4}
|g(t, x, \y(x), \nabla \y(x))| \le m \quad \text{ for }t\in[0,+\infty) \text{ and } \ x\in\Omega.
\end{equation}
Hence, for any $t\in[0,+\infty)$ and $\y\in X^\alpha$, we infer that
\begin{align*}
\|F(t,\y)\|^p & = \int_\Omega |g(t, x, \y(x),\nabla \y(x))|^p \d x \le m^p |\Omega|,
\end{align*}
and the inequality \eqref{sw} is satisfied with $K:= m |\Omega|^{1/p}$, which proves $(ii)$.
To verify that $F$ satisfies $(F1)$, take $t\in[0,+\infty)$ and $\y_1,\y_2\in X^\alpha$ and observe that assumption $(E1)$ implies
\begin{align*}
\|F(t, \y_1) - F(t, \y_2)\|^p & \le \int_\Omega |g(t, x, \y_1(x),\nabla \y_1(x)) - g(t, x, \y_2(x),\nabla \y_2(x))|^p \d x  \\
& \hspace{-10mm} \le L^p \left(\int_\Omega (|\y_1(x) - \y_2(x)| + |\nabla \y_1(x) - \nabla \y_2(x)|)^p \d x \right)\\
& \hspace{-10mm} \le 2^{p-1}L^p\left(\int_\Omega |\y_1(x) - \y_2(x)|^p \d x + \int_\Omega |\nabla \y_1(x)-\nabla \y_2(x)|^p \d x \right) \\
& \hspace{-10mm} \le 2^{p-1}L^p \|\y_1 - \y_2\|^p_{W^{1,p}(\Omega)} \le 2^{p-1}M^pL^p \|\y_1 - \y_2\|^p_\alpha.
\end{align*}
Consequently
\begin{align*}
\|F(t, \y_1) - F(t, \y_2)\| \le 2^{1-1/p}M L \|\y_1 - \y_2\|_\alpha \ffor t\in[0,+\infty), \ \y_1,\y_2\in X,
\end{align*}
and assumption $(F1)$ holds. We now check that $F$ is continuous. To this end let $(t_n)$ in $[0,+\infty)$ and $(\y_n)$ in $X^\alpha$ be a sequences such that $t_n \to t_0$ and $\y_n \to \y_0$ as $n \to \infty$. Suppose that $(n_k)$ is an increasing sequence of positive integers such that $n_k\to +\infty$ as $k\to +\infty$. In view of continuity of the inclusion $X^\alpha\subset W^{1,p}(\Omega)$, there is a subsequence $n_{k_l}$ of $(n_k)$ such that $\y_{n_{k_l}}(x) \to \y_0(x)$ and $\nabla \y_{n_{k_l}}(x) \to \nabla \y_0(x)$ as $l\to \infty$, for a.a. $x\in\Omega$. Then $$g(t_{n_{k_l}}, x, \y_{n_{k_l}}(x), \nabla \y_{n_{k_l}}(x)) \to g(t_0, x, \y_0(x), \nabla \y_0(x)) \aas l\to +\infty$$ for a.a. $x\in\Omega$. On the other hand, from the inequality \eqref{og4} it follows that
\begin{align*}
|g(t_{n_{k_l}}, x, \y_{n_{k_l}}(x), \nabla \y_n(x))| \le m \qquad\text{for a.a. } \ x\in\Omega \ \text{ and } \ l\ge 1.
\end{align*}
Therefore, by the dominated convergence theorem, we find that $$F(t_{n_{k_l}}, \y_{n_{k_l}}) \to F(t_0, \y_0) \aas l\to +\infty$$ in $X=L^p(\Omega)$, which completes the proof. \hfill $\square$

\subsection{Unique continuation property}

In this section we recall the facts concerning the unique continuation property. We start with the following definition.
\begin{definition}
We say that $\y\in W^{1,2}_{loc}(\Omega)$ is \emph{a distributional solution} of $\mathcal{A} \y = \lambda \y$, where $\lambda\in\R$, if
\begin{equation*}
\int_\Omega\sum_{i,j=1}^n a_{ij}(x)D_i \y(x) D_j\varphi(x) \d x = \int_\Omega\lambda \y(x)\varphi(x) \d x \ffor \varphi\in C^\infty_0(\Omega),
\end{equation*}
where $C^\infty_0(\Omega)$ is the set of smooth functions with compact support contained in $\Omega$.
\end{definition}

The following theorem is known as {\em the unique continuation property} for elliptic operators and is a consequence of  Theorem 1.1 from \cite{MR882069} and Proposition 3 from \cite{MR1151266}. For more detail see also \cite{MR1013816}, \cite{MR0599580}, \cite{MR1836814} and references contained therein.

\begin{theorem}\label{unif-cont}
Let $\lambda\in\R$ and let $\y\in W^{1,2}_{loc}(\Omega)$ be a distributional solution of the equation $\mathcal{A} \y = \lambda \y$ which is equal to zero on the set of positive measure. Then $\y(x) = 0$ for a.a. $x\in\Omega$.
\end{theorem}

\begin{corollary}\label{rem-jed-kon}
Assume that $\lambda = \lambda_k$, where $k\ge 1$, is an eigenvalue of $A_p$. If $\y\in \ker (\lambda I - A_p)\setminus \{0\}$ then the set $\{x\in\Omega \ | \ \y(x) = 0\}$ is of measure zero.
\end{corollary}
\noindent\textbf{Proof.} Assume that $\y\in D(A_p) \subset W^{2,p}(\Omega)$ satisfies $A_p \y = \lambda \y$, where $\lambda = \lambda_k$ for $k\ge 1$. Since $p \ge 2$ we see that $\y\in W^{1,2}_{loc}(\Omega)$. Furthermore for any $\varphi \in C^\infty_0(\Omega)$ we have
\begin{align*}
\int_\Omega\lambda \y(x)\varphi(x) \d x & = \int_\Omega \mathcal{A} \y (x) \varphi(x) \d x = - \int_\Omega \sum_{i,j=1}^n D_j(a_{ij}(x)D_i \y(x)) \varphi(x)\d x \\
& = \int_\Omega\sum_{i,j=1}^n a_{ij}(x)D_i \y(x) D_j\varphi(x) \d x,
\end{align*}
which proves that $\y$ is distributional solution of $\mathcal{A} \y = \lambda \y$. Since $\y \neq 0$, from Theorem \ref{unif-cont} it follows that the measure of the set $\{x\in\Omega \ | \ \y(x) = 0\}$ is equal to zero, which completes the proof. \hfill $\square$

\subsection{Resonant properties of Niemytzki operator}

In this section, our aim is to examine what assumptions should satisfy the mapping $f$ so that the associated Niemytzki operator $F$ meets the introduced earlier geometrical conditions. We start with the following theorem which says that well known \emph{Landesman-Lazer} conditions introduced in \cite{MR0267269} are actually particular case of conditions $(G1)$ and $(G2)$.

\begin{theorem}\label{lem-est2}
Suppose that the continuous functions $g_+,g_-\colon \Omega \to \mathbb{R}$ are such that
\begin{align*}
g_+(x) = \lim_{s \to +\infty} g(t,x,s,y) \quad\text{and}\quad g_-(x) = \lim_{s \to -\infty} g(t,x,s,y)
\end{align*}
for $x\in\Omega$, uniformly for $t\in[0,+\infty)$ and $y\in\R^n$. Let $B\subset X^\alpha_+\oplus X^\alpha_-$ be a subset bounded in the norm $\|\cdot\|_\alpha$. \\[5pt]
\makebox[5mm][l]{(i)} \parbox[t]{117mm}{Assume that 
$$\int_{\{\y>0\}} g_+(x) \y(x) \,d x  + \int_{\{\y<0\}} g_-(x) \y(x) \,d x > 0 \leqno{(LL1)}$$
for $\y\in\ker(\lambda I - A_p)\setminus\{0\}$. Then there is $R > 0$ such that for any $t\in[0,T]$ and $(\z, \y)\in B\times X_0$ with $\|\y\|_{L^2} \ge R$, we have the following inequality:
\begin{equation*}
\<F(t,\z + \y), \y\>_{L^2} > 0.
\end{equation*}}\\
\makebox[5mm][l]{(ii)} \parbox[t]{117mm}{Assume that 
$$\int_{\{\y>0\}} g_+(x) \y(x) \,d x + \int_{\{\y<0\}} g_-(x) \y(x) \,d x < 0 \leqno{(LL2)}$$
for $\y\in\ker(\lambda I - A_p)\setminus\{0\}$. Then there is $R > 0$ such that for any $t\in[0,T]$ and $(\z, \y)\in B\times X_0$ with $\|\y\|_{L^2} \ge R$, we have the following inequality:
\begin{equation*}
\<F(t, \z + \y), \y\>_{L^2} < 0.
\end{equation*}}
\end{theorem}
\noindent\textbf{Proof.} Since the proofs of points $(i)$ and $(ii)$ are analogous, we focus only on the first one. Suppose, contrary to the point $(i)$, that there are sequences $(t_n)$ in $[0,T]$, $(\z_n)$ in $B$  and $(\y_n)$ in $X_0$ such that $\|\y_n\|_{L^2} \to \infty$ when $n\to \infty$ and
\begin{equation}\label{g1b}
\<F(t_n, \z_n + \y_n), \y_n\>_{L^2} \le 0 \ffor n\ge 1.
\end{equation}
For $n\ge 1$, we define $\bar z_n:= \y_n/\|\y_n\|_{L^2}$. Since $X_0$ is finite dimensional space, with out loss of generality we can suppose that there is $\bar z_0\in X_0$ such that $\bar z_n \to \bar z_0$ in $L^2(\Omega)$ and $\bar z_n(x)\to \bar z_0(x)$ for a.a. $x\in\Omega$ as $n\to \infty$.
In view of the fact that $A_p$ has compact resolvents, Remark \ref{rem-pom2} $(d)$ says that $X^\alpha$ is compactly embedded in $X$. Therefore, the boundedness of $(\z_n)$ in $X^\alpha$, implies that this sequence is relatively compact in $X$. Hence, passing if necessary to a subsequence, we can also suppose that $\z_n \to \z_0$ in $X$ where $\z_0\in X = L^p(\Omega)$ and furthermore $\z_n(x) \to \z_0(x)$ for a.a. $x\in\Omega$ as $n\to \infty$. From \eqref{g1b}, we have
\begin{equation}\label{g7b}
\<F(t_n, \z_n + \y_n), \bar z_n - \bar z_0\>_{L^2} + \<F( t_n, \z_n + \y_n), \bar z_0\>_{L^2} \le 0
\end{equation}
for $n\ge 1$. Furthermore, by Lemma \ref{lem-nem-op} $(ii)$, the map $F$ is bounded, and hence the convergence $\bar z_n \to \bar z_0$ in $L^2(\Omega)$, implies that
\begin{equation}\label{g6b}
\<F(t_n, \z_n + \y_n), \bar z_n - \bar z_0\>_{L^2} \le \|F(t_n, \z_n + \y_n)\|_{L^2} \|\bar z_n - \bar z_0\|_{L^2} \to 0
\end{equation}
as $n\to +\infty$. If we define $\Omega_+:= \{x\in\Omega \ | \ \bar z_0(x) > 0\}$, $\Omega_-:= \{x\in\Omega \ | \ \bar z_0(x) < 0\}$ and $\bar c_n = \z_n + \y_n$, then
\begin{equation}
\hspace{-2mm}\begin{aligned}\label{g5b}
& \<F(t_n, \z_n + \y_n), \bar z_0\>_{L^2} = \int_\Omega g(t_n,x, \bar c_n(x),  \nabla \bar c_n(x) ) \bar z_0(x) \d x \\ & = \int_{\Omega_+} g(t_n,x, \bar c_n(x),  \nabla \bar c_n(x) ) \bar z_0(x) \d x \ + \int_{\Omega_-} g(t_n,x, \bar c_n(x),  \nabla \bar c_n(x) ) \bar z_0(x) \d x
\end{aligned}
\end{equation}
for $n\ge 1$. Observe that the equation
$$ \bar c_n(x) = \z_n(x) + \y_n(x) =  \z_n(x) + \|\y_n\|_{L^2} \bar z_n(x) \ \ \text{ for a.a. } \ x\in\Omega_+ \text{ and } n\ge 1$$
leads to the convergence
\begin{equation*}
\bar c_n(x) = \z_n(x) + \y_n(x) \to +\infty \ \ \text{ for a.a. } \ x\in\Omega_+ \ \text{ as } n\to \infty,
\end{equation*}
which together with assumption $(E2)$ and dominated convergence theorem gives
\begin{equation*}
\int_{\Omega_+} g(t_n,x, \z_n(x) + \y_n(x),  \nabla \z_n(x) + \nabla \y_n(x) ) \bar z_0(x) \d x \to
\int_{\Omega_+} g_+(x) \bar z_0(x) \d x
\end{equation*}
when $n\to +\infty$. Proceeding in the similar way, we infer that
\begin{equation*}
\int_{\Omega_-} g(t_n,x, \z_n(x) + \y_n(x),  \nabla \z_n(x) + \nabla \y_n(x) ) \bar z_0(x) \d x \to
\int_{\Omega_-} g_-(x) \bar z_0(x) \d x
\end{equation*}
when $n\to +\infty$. Hence, combining this with \eqref{g5b} yields
\begin{multline*}
\<F(t_n, \z_n + \y_n), \bar z_0\>_{L^2} \to
\int_{\Omega_+} g_+(x) \bar z_0(x) \d x + \int_{\Omega_-} g_-(x) \bar z_0(x) \d x \aas n\to \infty.
\end{multline*}
Therefore, letting $n\to \infty$ in \eqref{g7b} and using \eqref{g6b}, we infer that
\begin{equation}
\int_{\Omega_+} g_+(x) \bar z_0(x) \d x + \int_{\Omega_-} g_-(x) \bar z_0(x) \d x \le 0,
\end{equation}
which contradicts condition $(LL1)$, because $\|\bar z_0\|_{L^2} = 1$. Thus the proof of point $(i)$ is completed.  \hfill $\square$ \\

The following lemma proves that assumptions $(G1)$ and $(G2)$ are also consequences of \emph{the strong resonance conditions} from \cite{MR713209}, \cite{MR597281}, \cite{MR1055536}.

\begin{theorem}\label{lem-est3}
Assume that there is a continuous function $g_\infty \colon \Omega \to \mathbb{R}$ such that
\begin{equation*}
g_\infty(x)  = \lim_{|s| \to +\infty} g(t,x,s,y)\cdot s
\end{equation*}
for $x\in\Omega$, uniformly for $t\in[0,+\infty)$ and $y\in\R^n$. Let $B\subset X^\alpha_+\oplus X^\alpha_-$ be a set bounded in the norm $\|\cdot\|_\alpha$. \\[5pt]
\makebox[5mm][l]{(i)} \parbox[t]{117mm}{If the following condition is satisfied
\begin{equation*}\leqno{(SR1)}
\quad\left\{\begin{aligned}  
& \text{there exists a function } q\in L^1(\Omega) \text{ such that } \\
& g(t,x,s,y)\cdot s \ge q(x) \text{ for } (t,x,s,y)\in[0,+\infty)\times\Omega\times\R\times\R^n
\text{ and } \\ & \int_\Omega g_\infty(x)\d x > 0,
\end{aligned}\right.
\end{equation*}
then there is $R > 0$ with the property that for any $t\in[0,T]$ and $(\y,\y)\in B\times X_0$ with $\|\y\|_{L^2}\ge R$, one has
\begin{equation*}
\<F(t, \z + \y), \y\>_{L^2} > 0.
\end{equation*}}\\
\makebox[5mm][l]{(ii)} \parbox[t]{117mm}{If the following condition is satisfied
\begin{equation*}\leqno{(SR2)}
\quad\left\{\begin{aligned}
& \text{there exists a function } q\in L^1(\Omega) \text{ such that } \\ 
& g(t,x,s,y)\cdot s \le q(x) \text{ for } (t,x,s,y)\in[0,+\infty)\times\Omega\times\R\times\R^n\text{ and } \\ & \int_\Omega g_\infty(x)\d x < 0,
\end{aligned}\right.
\end{equation*} 
then there is $R > 0$ with the property that for any $t\in[0,T]$ and $(\z,\y)\in B\times X_0$ with $\|\y\|_{L^2}\ge R$, one has:
\begin{equation*}
\<F(t,\z + \y), \y\>_{L^2} < 0.
\end{equation*}}
\end{theorem}
  
\noindent\textbf{Proof.} It suffices to prove the first point, as the proof of the second one goes analogously. We argue by contradiction and assume that there are sequences $(t_n)$ in $[0,T]$, $(\z_n)$ in $B$ and $(\y_n)$ in $X_0$ such that $\|\y_n\|_{L^2} \to +\infty$ and
\begin{equation}\label{eq:39b}
\<F(t_n, \z_n + \y_n), \y_n\>_{L^2} \le 0 \ffor n\ge 1.
\end{equation}
Since $B\subset X^\alpha$ is a bounded set and the inclusion $X^\alpha\subset X$ is compact, passing if necessary to subsequence, we can assume that there is $\z_0\in X$ such that $\z_n \to \z_0$ in $X$ and $\z_n(x)\to \z_0(x)$ for a.a. $x\in\Omega$ as $n\to +\infty$. For any $n\ge 1$, define $\bar z_n:= \y_n/\|\y_n\|_{L^2}$. Since $X_0$ is a finite dimensional space we can also assume that there is $\bar z_0\in X_0$ such that $\bar z_n\to \bar z_0$ and $\bar z_n(x) \to \bar z_0(x)$ for a.a. $x\in\Omega$ as $n\to +\infty$. Put $\bar c_n:=\z_n + \y_n$ for $n\ge 1$ and take $x\in\Omega_+:=\{x\in\Omega \ | \ \bar z_0(x) > 0\}$. Then
\begin{equation}\label{eq:41b}
    \bar c_n(x) = \z_n(x) + \y_n(x) = \z_n(x) + \|\y_n\|_{L^2} \bar z_n(x) \to +\infty,
\end{equation}
when $n\to +\infty$. If we take $x\in\Omega_-:=\{x\in\Omega \ | \ \bar z_0(x) < 0\}$ we infer that
\begin{equation}\label{eq:42b}
    \bar c_n(x) = \z_n(x) + \y_n(x) = \z_n(x) + \|\y_n\|_{L^2} \bar z_n(x) \to -\infty
\end{equation}
when $n\to +\infty$. Using \eqref{eq:39b} we derive that
\begin{equation}\label{eq:40b}
\<F(t_n, \z_n + \y_n), \z_n + \y_n\>_{L^2} \le \<F(t_n, \z_n + \y_n), \z_n \>_{L^2}
\end{equation}
for any $n\ge 1$. Note that for the both conditions $(SR1)$ and $(SR2)$ we have
\begin{equation}
\begin{aligned}\label{rownww}
& \int_{\Omega_+} g(t_n,x,\bar c_n(x), \nabla \bar c_n(x) ) \bar c_n(x) \d x \ge -\|h\|_{L1} \ \text{ oraz } \\
& \int_{\Omega_-} g(t_n,x,\bar c_n(x),  \nabla \bar c_n(x) ) \bar c_n(x) \d x \ge -\|h\|_{L1} \ \text{ for } \ \ n\ge 1.
\end{aligned}
\end{equation}
Since $z_0 \neq 0$, from Corollary \ref{rem-jed-kon} it follows that the Lebesgue measure of the set $\Omega_0:=\{x\in\Omega \ | \ z_0(x) = 0\}$ is equal to zero. Therefore, applying the inequalities \eqref{rownww}, we infer that
\begin{equation*}
\begin{aligned}
& \liminf_{n\to +\infty}\<F(t_n, \z_n + \y_n),  \z_n + \y_n\>_{L^2} = \liminf_{n\to +\infty}  \int_\Omega g(t_n,x,\bar c_n(x), \nabla \bar c_n(x) ) \bar c_n(x) \d x \\
& \hspace{30mm}\ge \liminf_{n\to +\infty} \int_{\Omega_+} g(t_n,x,\bar c_n(x), \nabla \bar c_n(x) ) \bar c_n(x) \d x \\
& \hspace{30mm}\qquad + \liminf_{n\to +\infty} \int_{\Omega_-} g(t_n,x,\bar c_n(x),  \nabla \bar c_n(x) ) \bar c_n(x) \d x.
\end{aligned}
\end{equation*}
According to the assumption of lemma
$$g(t_n,x,\bar c_n(x), \nabla \bar c_n(x) ) \bar c_n(x) \ge h(x) \ffor n\ge 1, \ \text{ oraz p.w. } \ x\in\Omega,$$
and hence, combining \eqref{eq:41b}, \eqref{eq:42b} and Fatou lemma gives
\begin{equation*}
\begin{aligned}
& \liminf_{n\to +\infty}\<F(t_n, \z_n + \y_n),  \z_n + \y_n\>_{L^2} = \liminf_{n\to +\infty}  \int_\Omega g(t_n,x,\bar c_n(x), \nabla \bar c_n(x) ) \bar c_n(x) \d x \\
& \hspace{30mm}\ge \int_{\Omega_+} \liminf_{n\to +\infty} g(t_n,x,\bar c_n(x), \nabla \bar c_n(x) ) \bar c_n(x) \d x \\
& \hspace{30mm}\qquad + \int_{\Omega_-} \liminf_{n\to +\infty} g(t_n,x,\bar c_n(x),  \nabla \bar c_n(x) ) \bar c_n(x) \d x \\
& \hspace{30mm} = \int_{\Omega_+} g_\infty(x) \d x + \int_{\Omega_-} g_\infty(x) \d x =
\int_\Omega g_\infty(x) \d x,
\end{aligned}
\end{equation*}
which in turn, implies that
\begin{equation}\label{eq:44b}
\liminf_{n\to +\infty}\<F(t_n, \z_n + \y_n),  \z_n + \y_n\>_{L^2} \ge \int_\Omega g_\infty(x) \d x.
\end{equation}
Since $\Omega$ is a bounded set, the inclusion $X\subset L^2(\Omega)$ is continuous. Hence there is $M > 0$ such that $$\|\y\|_{L^2} \le M \|\y\|_\alpha \ffor \y\in X^\alpha.$$ From the boundedness of $B$, it follows that there is a constant $r < +\infty$ such that $r:=\sup\{\| \z_n \|_{L^2} \ | \ n\ge 1\}$. Then, for any $n\ge 1$,
\begin{equation}\label{eq:43b}
    \<F(t_n, \z_n + \y_n),  \z_n \>_{L^2} \le \|F(t_n, \z_n + \y_n)\|_{L^2} \| \z_n \|_{L^2} \le r \|F(t_n, \z_n + \y_n)\|_{L^2}.
\end{equation}
Note that, from the assumptions of lemma, we have
\begin{equation}\label{zbiez}
\lim_{|s| \to +\infty} g(t,x,s,y) = 0
\end{equation}
for $x\in\Omega$, uniformly for $t\in[0,+\infty)$ and $y\in\R^n$. Therefore, combining  \eqref{eq:41b}, \eqref{eq:42b} and \eqref{zbiez}, yields
\begin{equation*}
    g(t_n,x,\bar c_n(x), \nabla \bar c_n(x)) \to 0 \quad \text{ for a.a. } \ x\in\Omega_+\cup\Omega_-.
\end{equation*}
Since $\Omega_0$ is of Lebesgue measure zero, the boundedness of $f$ (assumption $(E2)$) and dominated convergence there imply that
\begin{equation*}
\begin{aligned}
    \|F(t_n, \z_n + \y_n)\|_{L^2}^2 & = \int_{\Omega_+} |g(t_n,x,\bar c_n(x), \nabla \bar c_n(x))|^2 \d x  \\
    & \qquad + \int_{\Omega_-} |g(t_n,x,\bar c_n(x), \nabla \bar c_n(x) )|^2 \d x \to 0,
\end{aligned}
\end{equation*}
when $n\to +\infty$. Hence the inequality \eqref{eq:43b} implies
\begin{equation*}
    \<F(t_n, \z_n + \bar c_n),  \z_n - \y_n\>_{L^2} \to 0 \aas n\to +\infty,
\end{equation*}
which along with \eqref{eq:40b} and \eqref{eq:44b}, leads to
\begin{equation}
0 \ge \liminf_{n\to +\infty} \<F(t_n, \z_n + \y_n),  \z_n + \y_n\>_{L^2} \ge \int_\Omega g_\infty(x) \d x.
\end{equation}
This inequality contradicts the condition $(SR1)$ and hence the proof of point $(i)$ is completed. \hfill $\square$

\subsection{Existence of periodic solutions}

In this section we intend to provide applications to study the existence of $T$-periodic solutions for particular differential equations being at resonance at infinity. It is worth nothing that the similar results were obtained by other authors, for example in \cite{MR0267269} \cite{MR0513090}, \cite{Aleks}, \cite{MR1766179}, \cite{MR1342038}. As a novelty we can recognize the fact that we examine equations were gradient is involved with the nonlinearity. To be more precise we shall consider parabolic equations of the form
\begin{equation}\label{A-eps-res-a}
v_t(t, x) = - \mathcal{A} \, v(t,x) + \lambda  v(t,x) + g(t,x, v(t,x), \nabla v(t,x)), \ \ t > 0, \ x\in\Omega
\end{equation}
where $\lambda\in\R$, and $f:[0,+\infty)\times\Omega\times\R\times\R^n\to\R$ is a continuous map satisfying assumptions $(E1)-(E3)$ and \\[5pt]
\noindent\makebox[22pt][l]{$(E4)$} \parbox[t][][t]{118mm}{there exists $T > 0$ such that $g(t,x,s,y) = g(t + T,x,s,y)$ for $t\in[0,+\infty)$, $x\in\Omega$, $s\in\R$, $y\in\mathbb{R}^n$.} \\[5pt] 
This equation may be written in the abstract form as
\begin{equation}\label{row-dr}
\dot v(t)  = - A_p v(t) + \lambda v(t) + F (t,v(t)), \qquad  t > 0.
\end{equation}
\begin{definition}
Let $J\subset \R$ be an interval. We say that $v:J\to X^\alpha$ is a solution of the equation \eqref{A-eps-res-a}, if $u$ is a mild solution of \eqref{row-dr}. 
\end{definition}
\noindent From Lemma \ref{lem-nem-op} it follows that $F$ satisfies $(F1)$ and $(F2)$. Hence Theorem \ref{th-exist1} implies that for any $\y_0\in X^\alpha$, equation \eqref{row-dr} admits a mild solution $u(\,\cdot\,;\y_0):[0,+\infty)\to X^\alpha$ such that $u(0;\y_0) = \y_0$. Define \emph{the translation along trajectories operator} $\Phi_T: X^\alpha \to X^\alpha$ associated with \eqref{row-dr} by $$\Phi_T(\y) := u(T; \y) \ffor \y\in X^\alpha.$$
Then $\Phi_T$ is a completely continuous map as a result of Theorems \ref{tw-con-21} and \ref{tw-exi-con-comp}.
We proceed to applications of the results obtained in previous sections to study the existence of periodic solutions. We start with the following \emph{criterion with Landesman-Lazer conditions}.
\begin{theorem}\label{crit-per-lan} 
Assume that there are continuous functions $g_\pm\colon \Omega \to \mathbb{R}$ such that
\begin{equation*}
g_+(x) = \lim_{s \to +\infty} g(t,x,s,y) \quad\text{and}\quad g_-(x) = \lim_{s \to -\infty} g(t,x,s,y)
\end{equation*} 
for $x\in\Omega$, uniformly for $t\in [0,+\infty)$ and $y\in\R^n$. If $\lambda=\lambda_k$ for some $k\ge 1$ and either $(LL1)$ or $(LL2)$ is satisfied, then the equation \eqref{A-eps-res-a} admits a $T$-periodic solution.
\end{theorem}

In the proof we use the following \emph{index formula with Landesman-Lazer conditions}, which is an immediate consequence of Remark \ref{rem-pom2} $(a)$ and Theorems \ref{lem-est2} and \ref{th-reso-m}.

\begin{theorem}\label{th-reso}
Under the assumptions of Theorem \ref{crit-per-lan} there is $R > 0$ such that $\Phi_T(\y)\neq \y$ for $\y\in X^\alpha$ with $\|\y\|_\alpha \ge R$ and 
\begin{enumerate}
\item[(i)] $\mathrm{deg_{LS}}(I - \Phi_T, B(0,R)) = (-1)^{d_k}$, if $(LL1)$ holds; \\[-9pt]
\item[(ii)] $\mathrm{deg_{LS}}(I - \Phi_T, B(0,R)) = (-1)^{d_{k-1}}$, if $(LL2)$ holds. 
\end{enumerate}
\end{theorem}

\noindent\textbf{Proof of Theorem \ref{crit-per-lan}.} By Theorem \ref{th-reso} and the existence property of topological degree, we see that each of the conditions $(LL1)$ or $(LL2)$ implies that there is $\y_0\in X^\alpha$ such that $\Phi_T(\y_0) = \y_0$. In view of assumption $(E4)$ we infer that $F(t,\y) = F(t+T, \y)$ for $t\ge 0$ and $\y\in X^\alpha$, which implies that $\y_0$ is a starting point of a $T$-periodic solution of \eqref{row-dr} being, by definition, a $T$-periodic solution of \eqref{A-eps-res-a}. \hfill $\square$ \\

Now we prove \emph{the criterion with strong resonance conditions}.

\begin{theorem}\label{crit-per-sr}
Let $\Omega\subset\R^n$ where $n\ge 3$, be an open bounded set and assume that there is a continuous function $g_\infty \colon \o\Omega \to \mathbb{R}$ such that
\begin{equation*}
g_\infty(x)  = \lim_{|s| \to +\infty} g(t,x,s,y)\cdot s
\end{equation*}
for $x\in\Omega$, uniformly for $t\in [0,+\infty)$ and $y\in\R^n$. If $\lambda=\lambda_k$ for some $k\ge 1$ and either condition $(SR1)$ or $(SR2)$ is satisfied, then the equation \eqref{A-eps-res-a} admits a $T$-periodic solution.
\end{theorem}

In the proof we use \emph{the index formula with strong resonance conditions}, which is a direct consequence of Remark \ref{rem-pom2} $(a)$ and Theorems \ref{lem-est3} and \ref{th-reso-m}. 

\begin{theorem}\label{th:2}
Under the assumption of Theorem \ref{crit-per-sr} there is $R > 0$ such that $\Phi_T(\y)\neq \y$ for $\y\in X^\alpha$ with $\|\y\|_\alpha \ge R$ and
\begin{enumerate}
\item[(i)] $\deg(I - \Phi_T, B(0,R)) = (-1)^{d_k}$, if $(SR1)$ holds; \\[-9pt]
\item[(ii)] $\deg(I - \Phi_T, B(0,R)) = (-1)^{d_{k-1}}$, if $(SR2)$ holds.
\end{enumerate}
\end{theorem}
\noindent\textbf{Proof of Theorem \ref{crit-per-sr}.} By Theorem \ref{th:2} and the existence property of topological degree, we see that each of the conditions $(SR1)$ and $(SR2)$ implies the existence of $\y_0\in X^\alpha$ such that $\Phi_T(\y_0) = \y_0$. In view of assumption $(E4)$ we deduce that $F(t,\y) = F(t+T, \y)$ for $t\ge 0$ and $\y\in X^\alpha$ and consequently $\y_0$ is a starting point of a $T$-periodic solution of \eqref{row-dr} and hence $T$-periodic solution of \eqref{A-eps-res-a}. \hfill $\square$ \\

\def\cprime{$'$} \def\polhk#1{\setbox0=\hbox{#1}{\ooalign{\hidewidth
  \lower1.5ex\hbox{`}\hidewidth\crcr\unhbox0}}} \def\cprime{$'$}
  \def\cprime{$'$} \def\cprime{$'$}
\providecommand{\bysame}{\leavevmode\hbox to3em{\hrulefill}\thinspace}
\providecommand{\MR}{\relax\ifhmode\unskip\space\fi MR }
\providecommand{\MRhref}[2]{%
  \href{http://www.ams.org/mathscinet-getitem?mr=#1}{#2}
}
\providecommand{\href}[2]{#2}

\parindent = 0 pt

\end{document}